\documentclass[review]{elsarticle}

\journal{Journal of \LaTeX\ Templates}

\usepackage{algorithm}
\usepackage{algorithmic}
\usepackage{color}
\usepackage{multirow}
\usepackage{supertabular}
\usepackage{longtable}
\usepackage{subfigure}
\usepackage{caption}
\usepackage{amsthm}
\usepackage{amsmath}
\usepackage{amssymb}
\usepackage{hyperref}
\usepackage{url}
\usepackage[normalem]{ulem}

\begin{document}

\begin{frontmatter}
\title{On the optimal pivot path of simplex method for linear programming based on reinforcement learning}


\author[mymainaddress]{Anqi Li}
\ead{lianqi20@mails.ucas.ac.cn}

\author[mymainaddress]{Tiande Guo}
\ead{tdguo@ucas.ac.cn}

\author[mymainaddress]{Congying Han\corref{mycorrespondingauthor}}
\cortext[mycorrespondingauthor]{Corresponding author}
\ead{hancy@ucas.ac.cn}

\author[mymainaddress]{Bonan Li}
\ead{libonan16@mails.ucas.ac.cn}

\author[mymainaddress]{Haoran Li}
\ead{lihaoran21@mails.ucas.ac.cn}

\address[mymainaddress]{School of Mathematical Sciences, University of Chinese Academy of Sciences, Shijingshan District, Beijing, China}


\begin{abstract}
Based on the existing pivot rules, the simplex method for linear programming is not polynomial in the worst case. Therefore the optimal pivot of the simplex method is crucial. This study proposes the optimal rule to find all shortest pivot paths of the simplex method for linear programming problems based on Monte Carlo tree search (MCTS). Specifically, we first propose the SimplexPseudoTree to transfer the simplex method into tree search mode while avoiding repeated basis variables. Secondly, we propose four reinforcement learning (RL) models with two actions and two rewards to make the Monte Carlo tree search suitable for the simplex method. Thirdly, we set a new action selection criterion to ameliorate the inaccurate evaluation in the initial exploration. It is proved that when the number of vertices in the feasible region is $C_n^m$, our method can generate all the shortest pivot paths, which is the polynomial of the number of variables. In addition, we experimentally validate that the proposed schedule can avoid unnecessary search and provide the optimal pivot path. Furthermore, this method can provide the best pivot labels for all kinds of supervised learning methods to solve linear programming problems.
\end{abstract}

\begin{keyword}
{simplex method}\sep {linear programming}\sep {pivot rules}\sep {reinforcement learning}
\MSC[2010] 00-01\sep  99-00
\end{keyword}

\end{frontmatter}


\section{Introduction}
The simplex method is a classical method for solving linear programming (LP) problems. Although it is a nonpolynomial time algorithm, its worst case rarely occurs and its average performance is better than that of polynomial time algorithms, such as the interior point method and ellipsoid method, especially for small-scale and medium-scale problems. Much research work has focused on making the simplex method a polynomial-time algorithm, but it has not been successful. The existing pivot rules can neither provide the optimal pivot paths for the simplex method nor make it a polynomial-time algorithm. In addition, the traditional design idea only applies to designing the pivot rule suitable for certain types of problems. There are no general ways to find the least number of pivot iterations for all types of linear programming. Our research goal is to design a general optimal pivot rule based on the inherent features of linear programming extracted by reinforcement learning (RL) that can be solved in polynomial time. This study is the first step toward achieving this goal. 

With the rise of machine learning (ML), ML-based technologies provides researchers with new ideas of pivot rules. Based on the deep Q-network (DQN)~\cite{DQN_1,DQN_2}, DeepSimplex~\cite{1} provides a pivot rule that can select the most suitable pivot rule for the current state between the Dantzig and steepest-edge rules. While another study~\cite{2} provides an instance-based method, the most suitable pivot rule for the current instance is learned among the five conventional pivot rules. The above two methods are based on several given pivot rules, and then learn the pivot rule scheduling scheme depending on the solution state or input instances. Therefore, the performance of these methods is heavily dependent on the supervised pivot rules. Unfortunately, owing to the lack of optimal labels, supervised pivot rules cannot extract the optimal pivot paths for the simplex method. 

In addition, the difficulty in determining the optimal pivot path lies in the information after several pivot iterations in the future. The existing solution state is insufficient for optimal future decisions. The most effective method is to appropriately assess the future situation before deciding to guide the best pivot. Fortunately, this idea is consistent with the Monte Carlo tree search (MCTS). Specifically, MCTS explores the trajectory in advance to evaluate and obtain future information to guide decision-making, significantly reducing the invalid search space and effectively guiding the best decision-making. Thus, the simplex method can effectively use future information to guide the current optimal pivot decision.

\begin{figure*}[t]
\centering
\includegraphics[width=1.00\columnwidth]{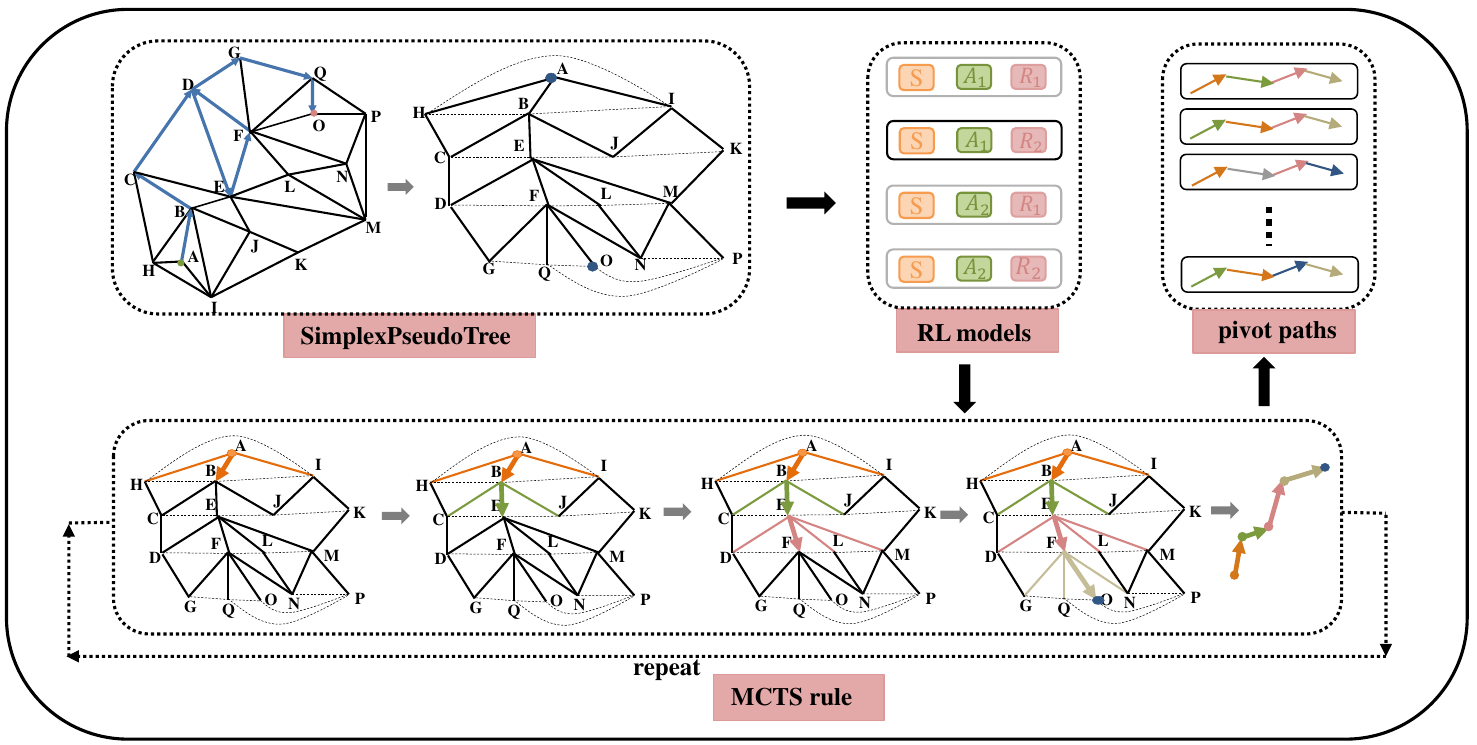}
\caption{Overview of the methodological framework in this paper. Firstly, we create SimplexPseudoTree to transform the simplex method applicable to reinforcement methods in Section 3. Next, four RL models are proposed in Section 4.1 based on the SimplexPseudoTree. Then we propose the MCTS rule to calculate all the shortest pivot paths in Section 4.2 and Section 4.3. Finally, we give thorough theory analysis for the MCTS rule in Section 5.}
\label{fig:framework_2}
\end{figure*}

Motivated by these observations, we propose to analyze and improve the simplex method in pivoting with the Monte Carlo tree search, further pushing forward the frontier of the simplex method for linear programming in a general way. This study focuses on four core aspects: (1) transforming the simplex method into a pseudo-tree structure, (2) constructing appropriate reinforcement learning models, (3) providing the MCTS rule to find all shortest pivot paths, and (4) giving thorough theory for the optimality and complexity of the MCTS rule, as shown in Figure~\ref{fig:framework_2}.

First, transforming the simplex method into the tree search mode is the premise for applying the Monte Carlo tree search method. Considering the connectivity and acyclicity characteristics, the tree structure can effectively avoid the generation of cycles in exploration paths. In this way, it ingeniously avoid repetition of the basis variables in exploration. To construct an imitative tree structure, SimplexPseudoTree, we propose to regard current states as nodes and the corresponding pivoting process as edges. The original linear programming instance is taken as the root node, and optimal solutions correspond to leaf nodes. Furthermore, we remove edges between nodes in the same layer and pivots from the deep layer to the shallow layer. Under this construction, the problem of finding the shortest path from the root node to the leaf node is equivalent to that of finding the optimal pivot rule.

Subsequently, we transform the problem of finding the optimal pivot path into four reinforcement-learning models. In particular, the following two novel action spaces and two reward functions are introduced. 
\begin{itemize}
\item Action set: (1) non-basis variables whose reduced costs are less than zero and (2) non-basis variables whose reduced costs are not equal to zero.
\item Reward functions: (1) opposite of pivot iterations and (2) linear decay weight estimation of the average variation in the objective value caused by a single pivot.
\end{itemize}
Based on the proposed actions and rewards, we design four models of different forms. Through a comprehensive comparison, it is find that the model comprising action set (1) and reward functions (2) achieves the highest efficiency with the least computational cost. Furthermore, we construct a novel action selection criterion for the simplex method to ameliorate inaccurate evaluation in the initial exploration. Subsequently, we present the MCTS rule based on the Monte Carlo tree search to determine the optimal pivot for current states. 

In addition, the optimal pivot that corresponds to the minimum pivot iterations is not necessarily unique. However, current research has not provided a way to find multiple optimal pivot paths. Unlike deterministic pivot rules, our MCTS rule exhibits certain randomness in the exploration stage. Specifically, the proposed exploration criterion can introduce a controllable scale factor based on the upper confidence bounds (UCB) method. Thus, additional randomness is added to the balance between the estimated value and explorations brought by the UCB algorithm. Additionally, the randomness of the MCTS rule can guide the selection of different actions to achieve the minimum pivot iterations under the guidance of optimality. 

Consequently, we prove the optimality and completeness of the MCTS rule. And we also prove the polynomial complexity of the optimal pivot iterations when the number of vertices in the feasible region is $C_n^m$. Concretely, the MCTS rule can find all the shortest pivot paths according to Wiener-khinchin law of large Numbers. Firstly, the MCTS rule can find the optimal pivot path when explorations approaches infinity. Then the MCTS rule can find all the different pivot paths when executions approach infinity. Additionally, from the perspective of combinatorial numbers, we prove that the minimum pivot iterations is polynomial of variables when the number of vertices in the feasible region is $C_n^m$. We also verify the polynomial iterations from the geometric perspective.

Given the above four aspects, we present a novel MCTS rule that provides all the shortest pivot paths. Additionally, we can label massive instances with little cost for the supervised pivot rule based on the proposed MCTS rule. Comprehensive experiments on the NETLIB benchmark and random instances demonstrated the efficacy of the MCTS rule. It is worth noting that compared with the minimum pivot iterations achieved by other popular pivot rules, our result is only 54.55\% for random instances and 49.06\% for NETLIB. 

Our main contributions are as follows: 
\begin{itemize}
\item Construct the SimplexPseudoTree to ensure that MCTS can be applied to the simplex method while avoiding duplicate bases.
\item Propose the MCTS rule to determine all the optimal pivot sequences.
\item Provide a method to obtain the optimal pivoting labels for the supervised pivot rule within the allowable range of the calculation cost.
\item Give comprehensive theory for the optimality and complexity of the MCTS rule.
\end{itemize}

The remainder of this paper is organized as follows. Sections 2 introduces the background and related works. Section 3 introduces the SimplexPseudoTree to translate the simplex method for applying RL methods. Section 4 presents the optimal MCTS rule for the simplex method. We prove the optimality and complexity of the proposed optimal pivot rule in Section 5. Section 6 presents experimental results. The conclusions and further works are presented in Sections 7 and 8, respectively.

\section{Background and Related Work}
\paragraph{\textbf{LP Problem}} Linear programming is a type of optimization problem where the objective function and constraints are linear. The standard form of the simplex method is as follows:
\begin{equation}\label{eq:LP_form_stand}
\begin{aligned}
&\min{c^{T}x} \\
&s.t.\ Ax=b \\
&\ \ \ \ \ x\geq0,
\end{aligned}
\end{equation}
where $c\in\mathbb{R}^{n}$ is the objective function coefficient, $A\in\mathbb{R}^{m\times n}$ is the constraint matrix, $b\in {\mathbb{R}}^{m}$ is the right-hand side, and all variables take continuous values in the feasible region. The purpose of linear programming is to find a solution that minimizes the value of the objective function in the feasible region, i.e., the so-called optimal solution, while the corresponding objective value is the optimal value. 

\paragraph{\textbf{Simplex Method}} The simplex method is clear and easy to understand. After providing an initial feasible solution, the pivoting process includes three parts: variable division, selection of the entering basis variable, and derivation of the leaving basis variable~\cite{DantzigRule}. The variable division step divides the variable $x\in\mathbb{R}^{n}$ into $x=\left[{x_B}^T,{x_N}^T\\ \right]^T$. Correspondingly, divide $A=[B,N]$ and $c=\left[{c_B}^T,{c_N}^T\right]^T$, and each column of the coefficient matrix $B$ corresponding to $x_B$ is required to be linearly independent. At this time, constraint $Ax=b$ can be written as
\begin{equation}\label{eq:LP_simplex_1}
\begin{aligned}
x_B=B^{-1}b-B^{-1}N x_{N}.
\end{aligned}
\end{equation}
We can obtain Formula (\ref{eq:LP_simplex_2}) by substituting the constraint into the objective function.
\begin{equation}\label{eq:LP_simplex_2}
\begin{aligned}
c^{T}x=c_{B}^{T}B^{-1}b+(c_{N}^{T}-c_{B}^{T}B^{-1}N)x_N
\end{aligned}
\end{equation}

As the first term is a constant value, the objective function value is determined only by the second term
\begin{equation}\label{eq:LP_simplex_3}
\begin{aligned}
\bar{c}^T=c_{N}^{T}-c_{B}^{T}B^{-1}N,
\end{aligned}
\end{equation}
which is called reduced costs. The components of $x_{N}$ corresponding to the part of Formula (\ref{eq:LP_simplex_4}) are non-basis variables that can cause a decrease in the objective function. 
\begin{equation}\label{eq:LP_simplex_4}
\begin{aligned}
J=\{j|\bar{c}_{j}<0\}
\end{aligned}
\end{equation}
Therefore, the selection of the entering basis variable is the process of selecting a basis variable from the non-basis variables mentioned above. Different pivot rules provide different methods for selecting the entering basis variable. In other words, the essence of the pivot rule of the simplex method is to convert a certain column between bases corresponding to feasible solutions. Accordingly, the feasible region polyhedron starts from the initial solution vertex, and each pivot corresponds to a step transition between adjacent vertices until it reaches the optimal solution vertex. When the basis variables are determined, we can use Formula (\ref{eq:LP_simplex_5}) to derive the leaving basis variable. 
\begin{equation}\label{eq:LP_simplex_5}
\begin{aligned}
x_{B}=B^{-1}b-B^{-1}Nx_{N} \geq 0
\end{aligned}
\end{equation}
After the initial feasible solution is given, the pivot process is repeated until the basis corresponding to the optimal solution is obtained, i.e., the end of the simplex method. 

\paragraph{\textbf{Classical Pivot Rules}} The pivot rule provides direction for the exchange of the basis variables of the simplex method. The simplex method has several classical pivot rules. The Danzig rule~\cite{DantzigRule} is to select the component corresponding to the most negative reduced cost as the entering basis variable, i.e., choose 
\begin{equation}\label{eq:LP_pivot_1}
\begin{aligned}
x_{B}=B^{-1}b-B^{-1}Nx_{N} \geq 0.
\end{aligned}
\end{equation}
The Bland rule~\cite{BlandRule} is to select the component with the smallest index from the variables with reduced costs less than zero, i.e., choose
\begin{equation}\label{eq:LP_pivot_2}
\begin{aligned}
\tilde{J} \in \mathop{\arg\min} \{j|\bar{c}_j<0 \}.
\end{aligned}
\end{equation}
The steepest-edge rule~\cite{steepestEdgeRule_1,steepestEdgeRule_2} uses the columns of the corresponding non-basis matrix and the basis matrix to standardize the reduced costs and selects the column corresponding to the smallest component, i.e., choose
\begin{equation}\label{eq:LP_pivot_3}
\begin{aligned}
\tilde{J} \in \mathop{\arg\min} \{\frac{\bar{c}_j}{\Vert B^{-1}N_{j} \Vert}|\bar{c}_j<0 \}.
\end{aligned}
\end{equation}
The idea of the greatest improvement rule~\cite{GreatestImprovementRule} is to take the product of the reduced costs and the maximum increment of each non-basis variable as the evaluation standard and select the non-basis variable with the minimum value as the entering basis variable. Finally, the devex rule~\cite{DevexRule_1,DevexRule_2}, an approximation of the steepest-edge rule, uses an approximate weight to replace the norm in the evaluation criterion of the steepest-edge rule. Considering the number of pivot iterations, different pivot rules apply to different problem types. However, there is no universal pivot rule that can determine minimum pivot iterations for general linear programming instances. In addition, there is an LP instance so that the corresponding simplex method is not a polynomial algorithm for any pivot rules given above.

\paragraph{\textbf{Pivot Rules Based on Machine Learning}} In recent years, the development of machine learning has provided new ideas for combinatorial optimization. Specifically, T. Guo and C. Han et al.~\cite{Optimize_7,Optimize_8} give overviews of solving combinatorial optimization problems. ML-based methods gradually emerge to solve combinatorial optimization problems, such as knapsack~\cite{Optimize_4}, TSP~\cite{Optimize_5,UseMCTS_TSP_1} and P-median problem~\cite{Optimize_6}. Simultaneously, there are many ML-based methods~\cite{Optimize_3,Optimize_0,Optimize_1,Optimize_2} involving continuous optimization problems. In terms of the linear programming problem, there are two methods for improving the pivot rules of the simplex method based on the machine learning method. DeepSimplex~\cite{1} uses the idea of Q-value iteration to learn the best scheduling scheme of the Dantzig rule and steepest-edge rule. Another study~\cite{2} used a boost tree and a neural network to learn an instance-based adaptive pivot rule selection strategy based on five classical pivot rules. However, the performance of these two supervised methods is severely limited by the provided pivot rules. In general, it is difficult to obtain the minimum pivot iterations based on supervised learning without effective labels.

\paragraph{\textbf{Combinatorial Optimization Methods Based on MCTS}} With the emergence of AlphaGo~\cite{AlphaGo} and AlphaGo Zero~\cite{AlphaGoZero}, reinforcement learning  represented by MCTS has been widely used in many classical problems~\cite{UseRL_1,UseRL_2,UseRL_3,UseRL_4}. Combinatorial optimization problems based on the MCTS can be solved in two ways. A classical idea is to design a MCTS-based framework for various types of combinatorial optimization problems~\cite{UseMCTS_manyCo_1,UseMCTS_manyCo_2}. Another idea is to design a MCTS-based algorithm to solve a specific combinatorial optimization problem, such as the traveling salesman problem~\cite{UseMCTS_TSP_1,UseMCTS_TSP_2,UseMCTS_TSP_3,UseMCTS_TSP_4} and Boolean satisfiability problem~\cite{UseMCTS_SAT_1,UseMCTS_SAT_2,UseMCTS_SAT_3,UseMCTS_SAT_4}. We adopt the latter idea to find multiple optimal pivot paths for the simplex method based on MCTS.  

\section{Constructed SimplexPseudoTree Model}
Reinforcement learning involves making the agent interact with the environment in a trial-and-error manner. The results of trial-and-error are fed back to the agent in the form of rewards to guide agent behavior and achieve the goal of maximizing the rewards. However, the current and past solution states are inadequate for determining the optimal pivot paths. Actually, the future information obtained by executing the pivot makes the difference. In this case, Monte Carlo tree search is more effective in finding the optimal pivot. Therefore, constructing an imitative tree-search model for the simplex method is our primary goal.

\begin{figure*}[hbtp]
\centering
\includegraphics[width=1.00\columnwidth]{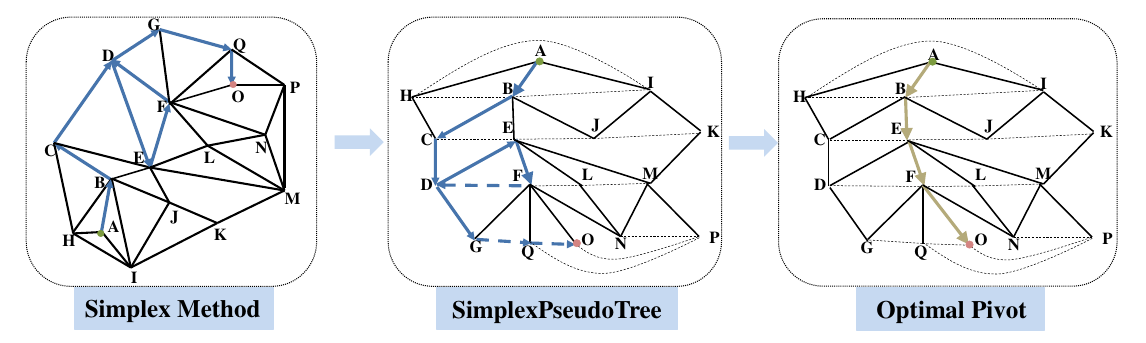}
\caption{Proposed SimplexPseudoTree for simplex method. The subgraph on the left shows the process of finding the optimal solution according to the simplex method. Middle subgraph is the SimplexPseudoTree corresponding to the instance. The subgraph on the right is the optimal pivot path found based on SimplexPseudoTree. }
\label{fig:Tree}
\end{figure*}

Duplicate basic variables lead to an unnecessary increase in pivot iterations in the simplex method. According to the reduced costs and pivot rules, the leaving basis variable cannot enter the basis into the immediate next pivot. Therefore, if the current solution state is considered a node and the pivot is considered an edge, the repetition of the basis implies that a circle appears in the exploration path, as shown in the left subgraph of Figure~\ref{fig:Tree}. Therefore, we must avoid base duplication in our pivot rule to minimize the number of pivot iterations consistent with the connectivity and acyclic properties of tree structures.

Inspired by the tree structure, SimplexPseudoTree is proposed to transform the simplex method for utilizing MCTS while avoiding the cycle during pivoting. The SimplexPseudoTree is constructed by considering the current states of the input instance as nodes and the feasible entering basis variables as edges. The initial linear programming instance is the root node, and the optimal solutions correspond to the leaf nodes. To minimize the number of pivot iterations, we only need to retain the pivot sequence for each node that first accesses it. Therefore, we remove the edges of same layers and pivots pointing to shallow layer nodes from deeper layers, as shown by the dotted line in Figure~\ref{fig:Tree}. In section 4.2, this idea is implemented by imposing significant negative rewards. It is worth noting that SimplexPseudoTree differs from the tree structure. SimplexPseudoTree allows multiple paths between two arbitrary nodes of SimplexPseudoTree because the pivot path is not unique. In general, the process of finding the optimal pivot rule is to find the shortest path between the root and leaf nodes for SimplexPseudoTree.

Using the idea of MCTS~\cite{AlphaGoZero,MCTS_1,MCTS_3}, we can evaluate the future situation for all candidate entering basis variables. Based on the explored information, we can significantly reduce the search space formed by all possible entering basis variables to find the optimal pivot paths by minimizing the number of pivot iterations. 

\section{Proposed RL Algorithm}
\subsection{RL Models of MCTS Rule}
Before applying reinforcement learning, this section presents the constructed state, action, and reward functions suitable for the simplex method. We provide two action-space definitions and two reward-function definitions. It is noteworthy that this is not a one-to-one correspondence. There are four combinations of RL models, as listed in Table~\ref{tab:model}.

\begin{table}[htpb]
\centering
\caption{Four reinforcement learning models constructed for the simplex method.}
\tiny
\resizebox{1\columnwidth}{!}{
\begin{tabular}{ccccc}
    \hline
    \textbf{Model} & \textbf{State} & \textbf{Action} & \textbf{Reward} \\
    \hline
    Model 1 & simplex tableaux & $A_{1}=\{i|\bar{c}_{i}<0\}$ & $R_{1}=-T$ \\
    Model 2 & simplex tableaux & $A_{2}=\{i|\bar{c}_{i}\neq 0\}$ & $R_{1}=-T$ \\
    Model 3 & simplex tableaux & $A_{1}=\{i|\bar{c}_{i}<0\}$ & $R_{2}=\frac{\sum_{i=1}^{N} w_{i}(c x_{i-1}-c x_{i})}{T}$ \\
    Model 4 & simplex tableaux & $A_{2}=\{i|\bar{c}_{i}\neq 0\}$ & $R_{2}=\frac{\sum_{i=1}^{N} w_{i}(c x_{i-1}-c x_{i})}{T}$ \\
    \hline
\end{tabular}}
\label{tab:model}
\end{table}

\paragraph{\textbf{State}}
We choose the simplex tableaux as the state representation of the reinforcement learning model. Simplex tableaux is the solution state representation of the simplex method on which traditional pivot rules rely. The tableaux contains $c$,$A$,$b$,$I_{B}$,$c_{B}$,$\bar{c}$, where $I_{B}$ are the column indexes corresponding to basis variables. The simplex tableaux completely represents the current solution state of the input instance and there is no redundant information.

\paragraph{\textbf{Two Action Sets}}
In the reinforcement learning model of the MCTS rule, actions corresponds to feasible entering basis variables. We have provided two definitions of action spaces. One action space contains variables with reduced costs of less than zero, i.e., non-basic variables whose objective value can be reduced by a one-step pivot, as shown in Formula (\ref{eq:action_leq}).
\begin{equation}
\label{eq:action_leq}
A_{1}=\{i|\bar{c}_{i}<0\}
\end{equation}
The other type of action space corresponds to variables whose reduced costs are not equal to zero. Compared with the former, this definition adds non-basic variables corresponding to reduced costs greater than zero (see Formula (\ref{eq:action_neq})). Although on the surface a variable with a reduced cost greater than zero does not make much sense. However, it represents the greedy idea of trying to find a path with fewer pivot iterations at the expense of the objective benefit in one step.
\begin{equation}
\label{eq:action_neq}
A_{2}=\{i|\bar{c}_{i}\neq 0\}
\end{equation}

\paragraph{\textbf{Two Reward Functions}}
We also provide two definitions for the reward function. The first reward function is defined as the opposite number of pivot iterations, which intuitively reflects the goal of minimizing the number of iterations, as shown in Formula (\ref{eq:reward_1}). One advantage of this is that in addition to having a minimum number of pivot iterations, the action selection guided by this reward is completely random, resulting in more randomness to find multiple pivot paths. 
\begin{equation}
\label{eq:reward_1}
R_{1}=-T
\end{equation}
The second reward function is defined by Formula (\ref{eq:reward_2}), where $T$ represents the maximum number of pivot iterations of the current episode, $i$ represents the $i^{th}$ pivot, $x_{i}$ is the locally feasible solution obtained from the $i^{th}$ pivot, and $w_{i} \in (0,1]$ is the weight. It is noteworthy that the proposed linear weight factor provides the weight of linear attenuation according to the depth from the root of the tree.
\begin{equation}
\label{eq:reward_2}
R_{2}=\frac{\sum_{i=1}^{N} w_{i}(c x_{i-1}-c x_{i})}{T} \ ,\ \ 
w_{i}=\frac{(T+1)-i}{T}
\end{equation}
Formula (\ref{eq:reward_2}) indicates a decrease in the linear weighted estimation of the objective value caused by a single pivot. In terms of minimizing pivot iterations, the two reward definitions are equivalent. However, as far as our problem is concerned, the second reward has the following two advantages: (1) compared with the first type of reward function, the dimensional feature of the change in the objective value is introduced; (2) the second reward function is more likely to choose the case in which the objective function changes significantly at the initial stage; therefore, even under the influence of the MCTS random exploration, this model can easily converge to the minimum pivot iterations.

\subsection{MCTS Rule}
\begin{figure*}[t]
\centering
\includegraphics[width=1.00\columnwidth]{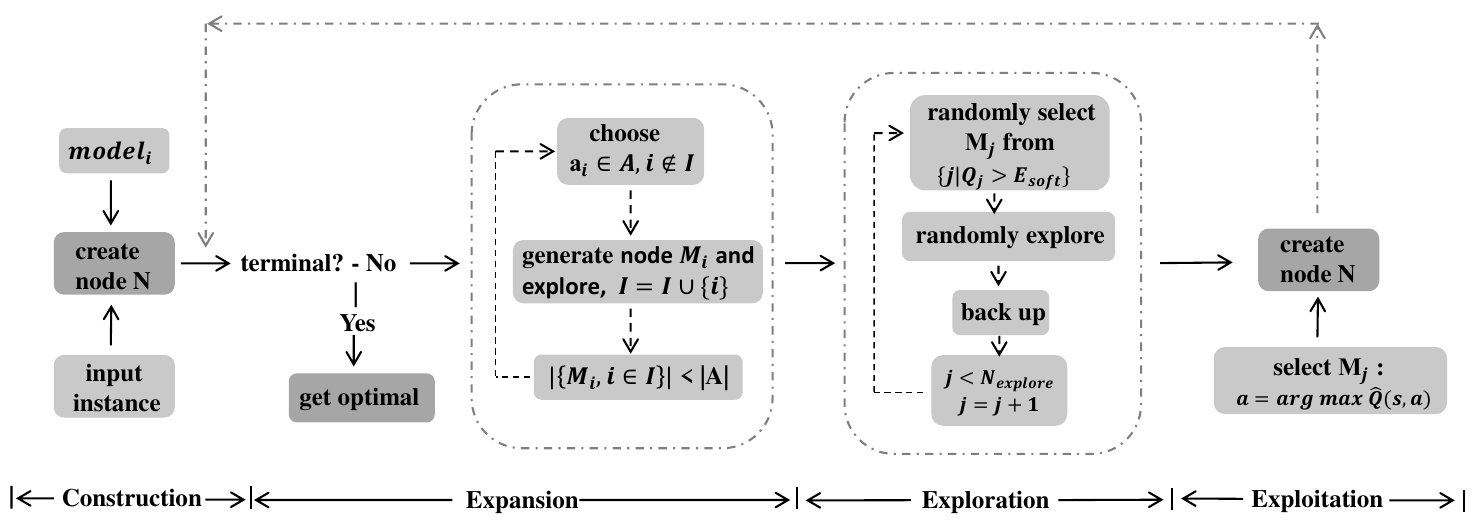} 
\caption{Algorithm flow diagram of the MCTS rule. }
\label{fig:flow_diagram}
\end{figure*}

Inspired by the idea of Monte Carlo tree search, we propose the MCTS rule, which can find the optimal pivot iterations for general linear programming instances. The entire process is divided into four stages: construction, expansion, exploration and exploration (Figure~\ref{fig:flow_diagram}). 

\paragraph{\textbf{Construction Stage}} 
First, we need to transform the simplex method into a structure with reinforcement learning model representation and an imitative tree-search pattern. The construction of SimplexPseudoTree is based on the schema in section 4.1, and the RL model is consistent with the strategy described in section 4.2. Specifically, the state of the present node represents the current solution stage, and each edge corresponds to an action in the action space. When an edge is selected from a node, the process enters the next solution stage via the corresponding pivot. Furthermore, a reward is obtained to evaluate the advantages and disadvantages of the current path while exploring the leaf nodes. Notably, we choose to enter the next stage unless the current node satisfies the optimality.

\paragraph{\textbf{Expansion Stage}} 
In the expansion stage, we randomly select actions in the action set without repetition to generate all child nodes of the selected node. The extraction process is shown in Formula (\ref{eq:equation_expansion_a}), where $A$ is the action set, $a_i$ is an action in the action set and $I$ is the set of currently selected action subscripts. 
\begin{equation}
\label{eq:equation_expansion_a}
randomly\ select\ from\ {\lbrace a_i\in A \vert i \notin I \rbrace}
\end{equation}
Considering the goal of minimizing the number of pivot iterations, the leaf node representing the optimality contains the  information we need. In terms of the optimal pivot path, the roles of the other nodes in the path are equivalent to those of the leaf nodes. In addition, the strategy evaluation method uses the empirical mean of the reward as the expectation of the reward. The state-action value function is defined as follows:
\begin{equation}
\label{eq:equation_3}
G_t=r_{T},
\end{equation}
\begin{equation}
\label{eq:equation_4}
Q_{\pi}(s,a)=\mathbb{E}_{\pi} [G_t \vert s_t=s,a_t=a], \  t\in \lbrace1,2,...,T\rbrace,
\end{equation} 
where $\left\{1,2,...,T \right\}$ denotes the set of visted subscripts in the current exploration path. After selecting all actions $\left\{a_i \in A\right\}$ and generating all possible child nodes $M=\left\{M_i,\  i\in I \right\}$, the process ends. 

\paragraph{\textbf{Exploration Stage}} 
The value function estimation is inaccurate when the number of explorations is low. At this time, the action corresponding to the maximum value function is not necessarily the real optimal choice, but slightly smaller values near it may do so. Therefore, we suggest relaxing the max operator in the initial stage to improve accuracy. The definitions are shown in Formulas (\ref{eq:equation_exploration_UCB}) and (\ref{eq:equation_exploration_Esoftmax}). The process of node selection is shown in Formula (\ref{eq:equation_exploration_select}), where the definition of $Q_i$ is consistent with that in the upper confidence bounds applied to the trees (UCT) algorithm~\cite{UCT}, and $v$ is the node explored currently. We randomly select actions from $ \lbrace a_i \vert Q_i>E_{soft} \rbrace$, then all actions in this episode are executed in a completely random manner.
\begin{equation}
\label{eq:equation_exploration_UCB}
Q_i=\frac{Q(v^{''})}{N(v^{''})}+C\sqrt{\frac{2\ln{N(v)}}{N(v^{''})}},\ v^{''}\in children\ of\ v
\end{equation}
\begin{equation}
\label{eq:equation_exploration_Esoftmax}
E_{soft}=\mathop{\min}\limits_{a_i \in A}{Q_i}+\alpha(\mathop{\max}\limits_{a_i \in A}{Q_i}-\mathop{\min}\limits_{a_i \in A}{Q_i})
\end{equation}
\begin{equation}
\label{eq:equation_exploration_select}
v^{'}=random \ select \ from \ \lbrace i \vert Q_i>E_{soft},a_i \in A \rbrace
\end{equation}
Considering that the number of pivot iterations required for subsequent duplicate nodes must be greater than the first one, this should be prohibited. Therefore, we only give real rewards to the nodes we encounter for the first time and punish repeated nodes by giving them huge negative rewards. During execution, when we encounter state $s$ and execute action $a$ for the first time in an episode, we add one to its count and increase the cumulative reward at that time, as shown in Formulas (\ref{eq:equation_6}) and (\ref{eq:equation_7}).
\begin{equation}
\label{eq:equation_6}
N(s,a) \leftarrow N(s,a)+1
\end{equation}
\begin{equation}
\label{eq:equation_7}
S(s,a) \leftarrow S(s,a)+G_t
\end{equation}
The state-action value function of the final state $s$ is in the form of Formula (\ref{eq:equation_8}). According to the law of large numbers, when the number of estimates tends to infinity, the value function tends to be close to that of the real strategy. When the number of iterations in this step reaches the preset threshold $N_{explore}$, the process terminates and enters the next stage.
\begin{equation}
\label{eq:equation_8}
Q(s,a)=S(s,a)/N(s,a)
\end{equation}

\paragraph{\textbf{Exploitation Phase}}
In the exploitation phase, we complete $N_{explore}$ exploration and estimate a reliable state-action value function. In this step, we select the action that maximizes the value function to generate nodes, as follows:
\begin{equation}
\label{eq:equation_exploitation_a}
a^{*}=\mathop{\arg\max}\limits_{a \in A}Q(s,a),
\end{equation}
\begin{equation}
\label{eq:equation_exploitation_v}
v^{*}=\mathop{\arg\max}\limits_{v^{''}\in children\ of\ v}\frac{Q(v^{''})}{N(v^{''})}.
\end{equation}
Subsequently, we must check whether the generated node has reached optimality; when it is not the optimal solution, we must return the node to the expansion stage and repeat the cycle.

\begin{algorithm}[H]
\renewcommand{\algorithmicrequire}{\textbf{Input:}}
\renewcommand{\algorithmicensure}{\textbf{Output:}}
\caption{the MCTS rule}
\label{alg:MCTS_rule}
\algsetup{linenosize=\tiny} \scriptsize

\begin{algorithmic}[1]
    \REQUIRE maximum number of explorations $N_{explore}$ and instance $I_{instance}$
	\ENSURE the optimal pivot path

    \STATE create node $N$, corresponding state $s_N$ and action set $A$ for $I_{instance}$
    \STATE $N_{(s_N,a_j)}=0,  S_{(s_{N},a_{j})}=0$ for $a_i \in A$
    \STATE $P=\emptyset $
    \WHILE{N is not terminal}
        \STATE $I=\emptyset $
        \WHILE{$\vert \lbrace M_i,i \in I \rbrace \vert<\vert A \vert$}
            \STATE choose $a_i \in A,i \notin I$
            \STATE generate node $M_i$ by executing $a_i$
            \STATE $I=I \bigcup \lbrace i \rbrace$
        \ENDWHILE 

        \STATE $j=0$    
        \WHILE{$j<N_{explore}$}
            \STATE randomly select $M_j$ from $\lbrace j \vert Q_j \geq E_{soft} \rbrace$
            \STATE randomly explore $M_j$ to get reward G
            \STATE $N_{(s_N,a_j)}=N_{(s_N,a_j)}+1$
            \STATE $S_{(s_{N},a_{j})}=S_{(s_{N},a_{j})}+G$
            \STATE $j=j+1$
        \ENDWHILE
    
        \STATE choose $\hat{a}=\mathop{\arg\max}\limits_{a\in A}{Q(s,a)}=\mathop{\arg\max}\limits_{a\in A}{S(s,a)/N(s,a)}$ 
        \STATE $P=P \cup \hat{a}$
        \STATE create node $\hat{N}$ by executing $\hat{a}$
        \STATE $N=\hat{N}$
        \STATE construct state $s_N$ and action set $A$ for node $N$
    \ENDWHILE
    \end{algorithmic}
\end{algorithm}

\subsection{Extracting Multiple Shortest Pivot Paths}
The optimal pivot paths correspond to different pivot sequences with minimum pivot iterations. Therefore, such an optimal path is not unique, which is important for the simplex method, but previous work is difficult to solve and ignores this point. Fortunately, the randomness of the MCTS rule is highly effective for finding multiple optimal paths. This randomness originates in the exploration stage. Specifically, the generation of the exploration trajectory depends on a random strategy. Thus, the estimated value brought about by limited exploration will be affected by randomness to a certain extent. Based on the randomness of the MCTS rule, our algorithm can select different actions that lead to optimal pivot paths in different execution processes. Therefore, multiple pivot sequences can be used to achieve optimization. Furthermore, we provide proof to ensure that each optimal pivot path can be found.

\section{Theoretical Analysis}
In this section, we conduct a detailed theoretical analysis of the MCTS rule from three perspectives. First, we prove that MCTS rule can make the optimal pivot decision at each step, so as to find the shortest pivot path. Then, considering the completeness of the algorithm, we can find all the optimal pivot paths when algorithm executions are sufficient. Finally, we prove that the pivot iterations under the MCTS rule is a polynomial of $n$ when the number of vertices in the feasible region is $C_n^m$.

\subsection{Optimality of the MCTS Rule}
We prove that the MCTS rule converges to the optimal pivot path when explorations approaches infinity. Considering that the expectation of the reward function in Models 1 and 2 will be affected by other episodes, it is not sufficient to reflect the real optimal pivot. Therefore, we first introduce the significance operator $Sig$ based on the idea of pooling in convolution, and then provide complete proof details.
\newtheorem{Definition}{Definition}[section]
\begin{Definition}[Rank Function]
    \label{def: rank func}
    Given a sequence of random variables $\mathrm{RS} := \{ X_1, X_2, \dots, X_n \} \subset \mathcal{X}$, it is sorted with an ascending order to get the order statistics sequence $\mathrm{RS_O} := \{ X_{(1)}, X_{(2)}, \dots, X_{(n)} \} $, where $X_{(i)}$ represents the $i$-th smallest random variable. Here we define the function $\operatorname{Rank}_{\mathrm{RS}}: \mathcal{X} \longrightarrow [n]$, where $[n] := \{ 1,2,\dots,n \}$, and $\operatorname{Rank}_{\mathrm{RS}}(X_i)$ is the index of $X_i$ in $\mathrm{RS_O}$.
\end{Definition}

\begin{Definition}[Significance Operator]
\label{Def1}
Given K groups of random variables sequences $\lbrace X^k_1,X^k_2,...,X^k_{n_k}\rbrace$, $k \in \lbrace 1,2,...,K \rbrace$, we define $\bar{X}^k$ and $X_{(n_k)}^k$ as the mean statistic and the maximum statistic of $k$-th sequence $\lbrace X^k_i\rbrace_{i=1}^{n_k}$, respectively:
$$\bar{X}^k = \frac{1}{n_k} \sum_{i=1}^{n_k} X_i^k,\ X_{(n_k)}^k = \max \{ X^k_1,X^k_2,...,X^k_{n_k} \},\ k \in \lbrace 1,2,...,K \rbrace.$$
For the mean statistics sequence $\mathrm{ES} := \{ \bar{X}^1,\bar{X}^2,\dots, \bar{X}^K \}$ and the maximum statistics sequence $\mathrm{MS} := \{ X_{(n_1)}^1,X_{(n_2)}^2,\dots, X_{(n_K)}^K \}$, we define the significance operator $\operatorname{Sig}: \mathcal{X} \longrightarrow \mathcal{X}$ as:
\begin{equation}
    \operatorname{Sig}(\bar{X}^k)=
    \begin{cases}
    \begin{array}{r}
        \bar{X}^k,\ \operatorname{Rank}_{\mathrm{ES}}(\bar{X}^k) = \operatorname{Rank}_{\mathrm{MS}}(X_{(n_k)}^k)\\
        X_{(n_k)}^k,\ \operatorname{Rank}_{\mathrm{ES}}(\bar{X}^k) \neq \operatorname{Rank}_{\mathrm{MS}}(X_{(n_k)}^k).
    \end{array}
    \end{cases}
\end{equation}
\end{Definition}

\newtheorem{thm}{\bf Theorem}[section]
\begin{thm}\label{thm1}
For Models 1 and 2, if we set $\alpha$ in Formula (\ref{eq:equation_exploration_Esoftmax}) as zero, the MCTS rule with significance operator $\operatorname{Sig}$ will converge to the optimal pivot rule with probability at least $1-\epsilon$, as long as
\begin{equation}
    N_{explore} \ge \frac{1}{\ln \left( 1 + \frac{1}{d_\mathcal{A}^*-1}\right)} \ln \left( \frac{1}{1-e^{\frac{1}{|P^*|} \ln\left( 1- \epsilon\right)} } \right) \approx O\left(  \ln \left( \frac{1}{\epsilon} \right) \right),
\end{equation}
where $P^*$ is the optimal pivot path and $|P^*|$ is its length. $d_\mathcal{A}^*:= \max_{s^*\in P^*} |\mathcal{A}_{s^*}|$ represents the dimension of the maximum action space along this path.
\end{thm}
\begin{proof}
For Model 1 and 2, the reward is defined as the opposite number of the pivot iterations. Therefore, pivot iterations corresponding to the maximum reward are equivalent to the shortest pivot path.

Without loss of generality, we consider the state node $s$ which represents a particular simplex tableaux. We denote the feasible action space in $s$ as $\mathcal{A}_s:=\{ a_1, a_2,\dots, a_{|\mathcal{A}_s|} \}$ where each action corresponds to a feasible pivot. The child nodes of $s$ are $\{ s_1, s_2, \dots, s_{|\mathcal{A}_s|}\}$ and the transition functions are $\mathbb{P}(s_i | s, a_j) = \mathbb{I}_{\{ i=j \}},\ \forall i,j\in \{1,2,\dots,|\mathcal{A}_s|\}$, where $\mathbb{I}$ is the indicator function. 

We denote $N_{explore}$ as the number of explorations from $s$. The random variable $N_i$ is the number of explorations from $s$ to $s_i$ and $\sum_{i=1}^{|\mathcal{A}_s|} N_i = N_{explore}$. 

If we set $\alpha$ in Formula (\ref{eq:equation_exploration_Esoftmax}) as zero, the set of selected actions is
\begin{equation}
\label{eq:equation_thm1_1}
 \left\{ i|Q(s,a_i) \ge E_{soft} \right\} = \left\{ i|Q(s,a_i) \ge \min_{a_j\in\mathcal{A}_s} Q(s,a_j)\right\} = \mathcal{A}_s,
\end{equation}
which indicates that all feasible pivot can be selected, i.e. $\mathbb{E}_{\pi}\left[N_i\right] > 0,\ \forall i\in \{1,2,\dots,|\mathcal{A}_s|\}$. In fact, our algorithm takes the exploration action by uniform policy from the set $\left\{ i|Q(s,a_i) \ge E_{soft} \right\}$, and then we have
\begin{equation}
\label{eq:equation_thm1_2}
\begin{aligned}
\mathbb{E}_{\pi}\left[N_i\right] &=  \sum_{i=1}^{|\mathcal{A}_s|} \mathbb{P}(s_i | s, a_j) \pi(a_j|s) N_{explore} \\
&= \pi(a_i|s) N_{explore} \\
&= \frac{N_{explore}}{|\left\{ i|Q(s,a_i) \ge E_{soft} \right\}|} \\
&= \frac{1}{|\mathcal{A}_s|}N_{explore} > 0, \qquad  \forall i\in \{1,2,\dots,|\mathcal{A}_s|\}.
\end{aligned}
\end{equation}
We can conclude that as long as $N_{explore} \ge \frac{1}{\ln \left( 1 + \frac{1}{|\mathcal{A}_s|-1}\right)} \ln \left( \frac{1}{\delta_1} \right)$, every child node can be accessed to with probability at least $1-\delta_1$, i.e. 
\begin{equation}\label{state visited}
    \begin{aligned}
        &\quad\ \mathbb{P}\left( \text{state $s_i$ visited} \right)\\
        &=\mathbb{P}{\left( N_i > 0 \right)}\\
        &= 1 - \mathbb{P}{\left( N_i = 0 \right)}\\
        &= 1 - \left( \frac{|\mathcal{A}_s|-1}{|\mathcal{A}_s|} \right)^{N_{explore}} \\
        &\ge 1-\delta_1.
    \end{aligned}
\end{equation}

We denote random variable $R_j^{a_i}$ as the reward of taking action $a_i$ for the $j^{th}$ time. Given $|\mathcal{A}_s|$ groups of independent identically distributed random variables sequences $\left\{ R_j^{a_i} \right\}_{j=1}^{N_i},\ i\in \{1,2,\dots,|\mathcal{A}_s|\}$, we apply the significance operator $\operatorname{Sig}$ in Definition \ref{Def1} and we have $\left\{ \operatorname{Sig}(\bar{R}^{a_i})\right\}_{i=1}^{|\mathcal{A}_s|}$. Note that the reward random variables $\left\{ R_j^{a_i} \right\}_{j=1}^{N_i}$ are independent and identically distributed, and we can apply Wiener-khinchin Law of Large Numbers: 
\begin{equation}
    \label{eq: LLN}
    \lim_{N_{explore} \to \infty} \bar{R}^{a_i} = \lim_{N_{i} \to \infty} \bar{R}^{a_i} =\lim_{N_{i} \to \infty} \frac{1}{N_i} \sum_{j=1}^{N_i} {R}^{a_i}_j = \mathbb{E}\left[ {R}^{a_i}_j \right],\quad a.s.
\end{equation}
Note that $Q(s,a_i)$ is the same as $\bar{R}^{a_i}$ with respect to the definition of $Q$.

For any pivot path $P = \{s^P_0, s^P_1, \dots, s^P_{|P|-1}\}$, we denote $d_\mathcal{A}^P$ as the dimension of the maximum action space along this path, i.e. $d_\mathcal{A}^P:= \max_{s^P\in P} |\mathcal{A}_{s^P}|$. 
By recursion, the MCTS rule can explore pivot path P with probability at least $1-\delta_2$ when $N_{explore} \ge \frac{1}{\ln \left( 1 + \frac{1}{d_\mathcal{A}^P-1}\right)} \ln \left( \frac{1}{1-e^{\frac{1}{|P|} \ln\left( 1- \delta_2\right)} } \right) \approx O\left( \frac{1}{\ln \left( 1 + \frac{1}{d_\mathcal{A}^P-1}\right)}  \ln \left( \frac{|P|}{\delta_2} \right) \right)$ , i.e. 

\begin{equation} \label{eq: access to each path}
    \begin{aligned}
        &\quad\ \mathbb{P}\left( \text{explore pivot path $P$} \right)\\
        &= \mathbb{P} \left( \bigcap_{s^P\in P} \left\{ N_{s^P} >0 \right\} \right)\\
        &= \prod_{s^P\in P} \mathbb{P} \left( N_{s^P} >0 \right)\\
        &= \prod_{s^P\in P} 1 - \left( \frac{|\mathcal{A}_{s^P}|-1}{|\mathcal{A}_{s^P}|} \right)^{N_{explore}} \\
        &\ge \left( 1 - \left( \frac{d_\mathcal{A}^P-1}{d_\mathcal{A}^P} \right)^{N_{explore}} \right)^{|P|}\\
        &\ge 1-\delta_2.
    \end{aligned}
\end{equation}
We have shown that our algorithm can explore each path $P$ when $N_{explore}$ is sufficiently large. As a result, we have
\begin{equation}
    \label{eq: maximum converge}
    \lim_{N_{explore} \to \infty} {R}^{a_i}_{(N_i)} = \lim_{N_{i} \to \infty} {R}^{a_i}_{(N_i)} =\lim_{N_{i} \to \infty} \max\left\{{R}^{a_i}_1, {R}^{a_i}_2,\dots, {R}^{a_i}_{N_i}\right\} = {R}^{a_i}_*,\quad a.s.
\end{equation}
where ${R}^{a_i}_*$ is the maximum reward which can be attained by taking action $a_i$. Denote $\mathrm{ES} := \{ \bar{R}^{a_1},\bar{R}^{a_2},\dots, \bar{R}^{a_{|\mathcal{A}|}} \}$ and the maximum statistics sequence $\mathrm{MS} := \{ R_{(N_1)}^{a_1},R_{(N_2)}^{a_2},\dots, R_{(N_{|\mathcal{A}|})}^{a_{|\mathcal{A}|}} \}$. Then, we have
\begin{equation}\label{eq: limit of sig}
    \lim_{N_{explore}\to\infty} \operatorname{Sig}(\bar{R}^{a_i})=
    \begin{cases}
    \begin{array}{l}
        \mathbb{E}\left[ {R}^{a_i}_j \right],\ \operatorname{Rank}_{\mathrm{ES}}(\bar{R}^{a_i}) = \operatorname{Rank}_{\mathrm{MS}}({R}^{a_i}_{(N_i)})\\
        {R}^{a_i}_*,\ \operatorname{Rank}_{\mathrm{ES}}(\bar{R}^{a_i}) \neq \operatorname{Rank}_{\mathrm{MS}}({R}^{a_i}_{(N_i)}),\quad a.s.
    \end{array}
    \end{cases}
\end{equation}

We define the mapping $\operatorname{Proj}: ES\cup MS \longrightarrow MS$, where $\operatorname{Proj}(\bar{R}^{a_i}) = {R}^{a_i}_{(N_i)},\ \forall \bar{R}^{a_i}\in ES$ and $\operatorname{Proj}({R}^{a_i}_{(N_i)}) = {R}^{a_i}_{(N_i)},\ \forall {R}^{a_i}_{(N_i)} \in MS$. Combined with Formula (\ref{eq: limit of sig}), we have
\begin{equation}
    \label{eq: limit of proj}
    \lim_{N_{explore}\to\infty} \operatorname{Proj}\circ\operatorname{Sig}\left(\bar{R}^{a_i} \right) = {R}^{a_i}_*,\quad a.s. \quad \forall i\in \{1,2,\dots,|\mathcal{A}_s|\}.
\end{equation}
We take action $\hat{a} \in \arg \max_{a_i\in\mathcal{A}_s} \operatorname{Proj}\circ\operatorname{Sig}\left(\bar{R}^{a_i} \right)$. According to the definition of $\operatorname{Proj}$ and $\operatorname{Sig}$, we can access the child node that attains the optimal reward to execute the next iteration. We have proven that the MCTS rule can find the optimal action from a starting state $s$. In the following, we will analysis the complexity to extract the whole optimal pivot path.

We assume the optimal pivot path is $P^* = \{s^*_0, s^*_1, \dots, s^*_{|P^*|-1}\}$ and denote $d_\mathcal{A}^*$ as the dimension of the maximum action space along this path, i.e. $d_\mathcal{A}^*:= \max_{s^*\in P^*} |\mathcal{A}_{s^*}|$.

By recursion, the MCTS rule will converge to the optimal pivot rule with probability at least $1-\epsilon$ when $N_{explore} \ge \frac{1}{\ln \left( 1 + \frac{1}{d_\mathcal{A}^*-1}\right)} \ln \left( \frac{1}{1-e^{\frac{1}{|P^*|} \ln\left( 1- \epsilon\right)} } \right) \approx O\left( \frac{1}{\ln \left( 1 + \frac{1}{d_\mathcal{A}^*-1}\right)}  \ln \left( \frac{|P^*|}{\epsilon} \right) \right)$, i.e. 
\begin{equation}
    \begin{aligned}
        &\quad\ \mathbb{P}\left( \text{find optimal pivot path $P^*$} \right)\\
        &= \mathbb{P} \left( \bigcap_{s^*\in P^*} \left\{ N_{s^*} >0 \right\} \right)\\
        &= \prod_{s^*\in P^*} \mathbb{P} \left( N_{s^*} >0 \right)\\
        &= \prod_{s^*\in P^*} 1 - \left( \frac{|\mathcal{A}_s|-1}{|\mathcal{A}_s|} \right)^{N_{explore}} \\
        &\ge \left( 1 - \left( \frac{d_\mathcal{A}-1}{d_\mathcal{A}} \right)^{N_{explore}} \right)^{|P^*|}\\
        &\ge 1-\epsilon.
    \end{aligned}
\end{equation}
Note that $d_\mathcal{A}$ does not exceed the number of feasible entering basis variables, i.e. $d_\mathcal{A} \le n-m$ and we also prove that $|P^*|$ is at most $\min(m,n-m)$ in Theorem \ref{thm3}, which indicates that $d_\mathcal{A}$ and $|P^*|$ are both linear according to the dimension of input. Therefore, we have that $N_{explore} \approx O\left( \frac{1}{\ln \left( 1 + \frac{1}{d_\mathcal{A}-1}\right)}  \ln \left( \frac{|P^*|}{\epsilon} \right) \right) \approx O\left( \ln \left( \frac{1}{\epsilon} \right) \right)$.

In summary, the MCTS rule will converge to the optimal pivot rule with the highest reward in model 1 and 2 when $N_{explore}$ is sufficiently large. 
\end{proof}

The above theorem concludes that the MCTS rule with significance operator converges to optimality under the condition of $\alpha = 0$. Additionally, the proof of convergence to the optimal strategy is given in the literature~\cite{UCT} for the case of $\alpha = 1$. It is also feasible to substitute them into the framework of our proof based on the upper and lower bounds of explorations of each action given in literature~\cite{UCT}. For other $\alpha$ values, we conducted an ablation study in the experimental section for verification. 

\subsection{Completeness of Multiple Pivot Paths}
The proposed MCTS rule is a random algorithm. Multiple executions may find different optimal pivot paths. Theorem \ref{thm2} gives a theoretical guarantee based on the Wiener-khinchin law of large Numbers. It is proved that the MCTS rule can find all the optimal pivot paths when executions are sufficient.
\begin{thm}\label{thm2}
For Model 1 and Model 2, if we set $\alpha$ in Formula (\ref{eq:equation_exploration_Esoftmax}) as zero, the MCTS rule with significance operator $\operatorname{Sig}$ can find all the optimal pivot paths provided that algorithm execution times $N_{exe}$ is sufficiently large.
\end{thm}
\begin{proof}
We use the same notations as ones in the proof of Theorem \ref{thm1}.
Set all optimal actions of the current pivot as $\mathcal{A}^*=\lbrace a^*_1,a^*_2,...,a^*_{|\mathcal{A}^*|}\rbrace$. When $|\mathcal{A}^*|=1$, this theorem holds by Theorem \ref{thm1}. We consider $|\mathcal{A}^*|\ge 2$ in the following proof. According to Formula (\ref{eq:equation_thm1_2}) 
\begin{equation}
\label{eq:equation_thm2_1}
\lim_{N_{explore} \to \infty} N_{a^*}=\lim_{N_{explore} \to \infty} \frac{1}{|\mathcal{A}|} N_{explore}=\infty, \ \forall a^* \in \mathcal{A}^*,
\end{equation}
where $\mathcal{A}$ represents the total action space, $N_{a^*}$ represents the number of explorations of the optimal action $a^*$.

According to Formula (\ref{eq: limit of proj}), we have 
\begin{equation}
    \label{eq: limit of opt}
    \lim_{N_{explore}\to\infty} \operatorname{Proj}\circ\operatorname{Sig}\left(\bar{R}^{a^*} \right) = {R}^{a^*}_* = \max_{i,j} R^{a_i}_j \quad a.s. \quad \forall a^* \in \mathcal{A}^*.
\end{equation}
Therefore, we select action from $\mathcal{A}^*$ by uniform policy, i.e. 
$\mathbb{P}({a^*_i})=\frac{1}{|\mathcal{A}^*|}.$ Then, we have that when $N_{exe} \ge \frac{1}{\ln \left( \frac{|\mathcal{A}^*|}{|\mathcal{A}^*|-1} \right)} \ln \left( \frac{|\mathcal{A}^*|}{\epsilon} \right)$, the MCTS rule can find all $a^*\in \mathcal{A}^*$, i.e.
\begin{equation}
    \label{eq:execute number inf}
    \begin{aligned}
        & \mathbb{P}\left( \text{find all } a^*\in \mathcal{A}^* \right)\\
        =& 1-\mathbb{P}\left( \bigcup_{i=1}^{|\mathcal{A}^*|} \left\{ \text{not find } a^*_i \right\} \right)\\
        \ge& 1-\sum_{i=1}^{|\mathcal{A}^*|} \mathbb{P}\left(\text{not find } a^*_i \right)\\
        =&1- |\mathcal{A}^*| \left( \frac{|\mathcal{A}^*| -1}{|\mathcal{A}^*| } \right)^{N_{exe}}\\
        \ge&1-\epsilon.
    \end{aligned}
\end{equation}
It indicates that when the number of algorithm executions $N_{exe}$ approaches infinity, each $a^*\in \mathcal{A}^*$ can be found. Then repeating the above process, the MCTS rule can find all the optimal pivot paths.
\end{proof}

\subsection{Complexity of the Optimal Pivot}
Theorem \ref{thm3} proves that the MCTS rule can find polynomial pivot iterations when the number of vertices in the feasible region is $C_n^m$.
\begin{thm}\label{thm3}
For the standard form of simplex method for linear programming as Formula (\ref{eq:LP_form_stand}), $P=\left\{ x \in\mathbb{R}^{n}| Ax = b, x\geq0 \right\}$ represents the feasible region, where $A\in\mathbb{R}^{m\times n}$ and $rank(A)=m$. Suppose the number of feasible vertices of $P$ is $C_n^m$. Then the shortest distance (the minimum hops) between any two feasible vertices of $P$ is $\min(m,n-m)$.
\end{thm}
\begin{proof}
For the convenience of proving, we first convert the simplex of the feasible region into the SimplexPseudoTree structure proposed in Section 3. In this way, each vertex only appears once in different layers, so the number of vertices in each layer adds up to the total number of vertices $N$. Compared with the root of SimplexPseudoTree, the nodes of layer $i$ is obtained from the root through the $i$ pivot iterations, e.i. $C_m^i C_{n-m}^i$. There we have 
\begin{equation}
\label{eq:equation_thm3_1}
\sum_{i=0}^{T}{C_m^i C_{n-m}^i}=N=C_n^m,
\end{equation}
where $T$ is the layers of the SimplexPseudoTree and also the longest distance to the root. We discuss it in two cases, based on the Vandermonde's identity of combination number 
\begin{equation}
\label{eq:equation_thm3_2}
\sum_{i=0}^{l}{C_a^i C_{b}^{l-i}}=C_{a+b}^{l}.
\end{equation}
When $n-m \leq m$, we have 
\begin{equation}
\label{eq:equation_thm3_3}
\sum_{i=0}^{k}{C_{m}^{i} C_{n-m}^{k-i}}=C_{n}^{k}.
\end{equation}
Let $k=n-m$, then 
\begin{equation}
\label{eq:equation_thm3_4}
\sum_{i=0}^{n-m}{C_{m}^{i}C_{n-m}^{n-m-i}}=
\sum_{i=0}^{n-m}{C_{m}^{i}C_{n-m}^{i}}=
C_{n}^{n-m}=C_{n}^{m},
\end{equation}
i.e. $T=n-m$. When $m < n-m$, we have 
\begin{equation}
\label{eq:equation_thm3_5}
\sum_{i=0}^{k}{C_{m}^{k-i} C_{n-m}^{i}}=C_{n}^{k}.
\end{equation}
Let $k=m$, then 
\begin{equation}
\label{eq:equation_thm3_6}
\sum_{i=0}^{m}{C_{m}^{m-i} C_{n-m}^{i}}=\sum_{i=0}^{m}{C_{m}^{i} C_{n-m}^{i}}=C_{n}^{m},
\end{equation}
i.e. $T=m$. In this way, $T=\min(m,n-m)$, which is a polynomial of $n$. Additionally, $T$ can represent the maximum hops of all shortest paths between any two vertices in the feasible region. Therefore, the shortest pivot iterations starting from any initial point in the feasible region is the polynomial of the number of variables when the number of vertices in the feasible region is $C_n^m$.
\end{proof}

From the perspective of combination number, Theorem \ref{thm3} proves that the shortest distance between any two vertices of the feasible region is $\min(m,n-m)$ when the number of vertices in the feasible region is $C_n^m$. This conclusion is also meaningful from the geometric perspective of the simplex method. Specifically, the number of vertices in the feasible region is $C_n^m$ means that any $m$ columns of the constraint matrix corresponds to a basis matrix of the simplex method. Therefore, for any two basis matrices $B_1$ and $B_2$, there can be at most $\min (m, n-m)$ different columns. We only need to exchange different columns respectively. In this way, the initial vertex $B_1$ is converted to $B_2$ through $min (m, n-m)$ pivot iterations. In other words, the optimal pivot iterations needed for any linear objective functions is $min (m, n-m)$ at most under the conditions of the above theorems.

We reveal that the optimal pivot of simplex method is polynomial when the number of vertices in the feasible region is $C_n^m$ from the perspective of theory and geometry respectively. Furthermore, the proposed MCTS rule can find the polynomial pivot iterations when the number of vertices in the feasible region is $C_n^m$.

\newtheorem{corollary}{Corollary}[section]
\begin{corollary}\label{cor1}
When the number of vertices in the feasible region is $C_n^m$, the MCTS rule ensures that the number of pivot iterations becomes the polynomial of the number of variables.
\end{corollary}

\section{Experiment}
\subsection{Datasets and Experiment Setting}
\paragraph{Datasets} We conduct experiments using the NETLIB$ \footnote{http://www.netlib.org/lp/data/index.html} $ benchmark~\cite{NetlibArtical} and random instances. The method for generating random instances is as follows. First, we write the equivalent form of Formula (\ref{eq:LP_form_stand}), 
\begin{equation}\label{eq:LP_form_equal}
\begin{aligned}
&\min c^T x+\textit{0}^T x_0 \\
&s.t.\ Ax+I x_0=b \\
&\ \ \ \ \ x\geq0, x_0\geq0,
\end{aligned}
\end{equation}
where $x_0\in\mathbb{R}^{m}$, $\textit{0}\in\mathbb{R}^{m}$ is a zero vector and $I\in\mathbb{R}^{m\times m}$ is an identity matrix. Thus, we provide a setting for the initial feasible basis $B_0=I$, where $n$ basis variables correspond to $n$ columns of constraint matrix $A$ at the right end. Each component of $A,b$ and $c$ is a random number of uniform distributions of $[0,1000)$. For the constraint matrices of non-square matrices, the rows and columns are obtained from a uniform distribution of [0,800). 

\begin{figure}[htpb]
    \centering
    \subfigure[\label{fig:a}rand 50 $\times$ 50]
    {\includegraphics[width=1.00\columnwidth]{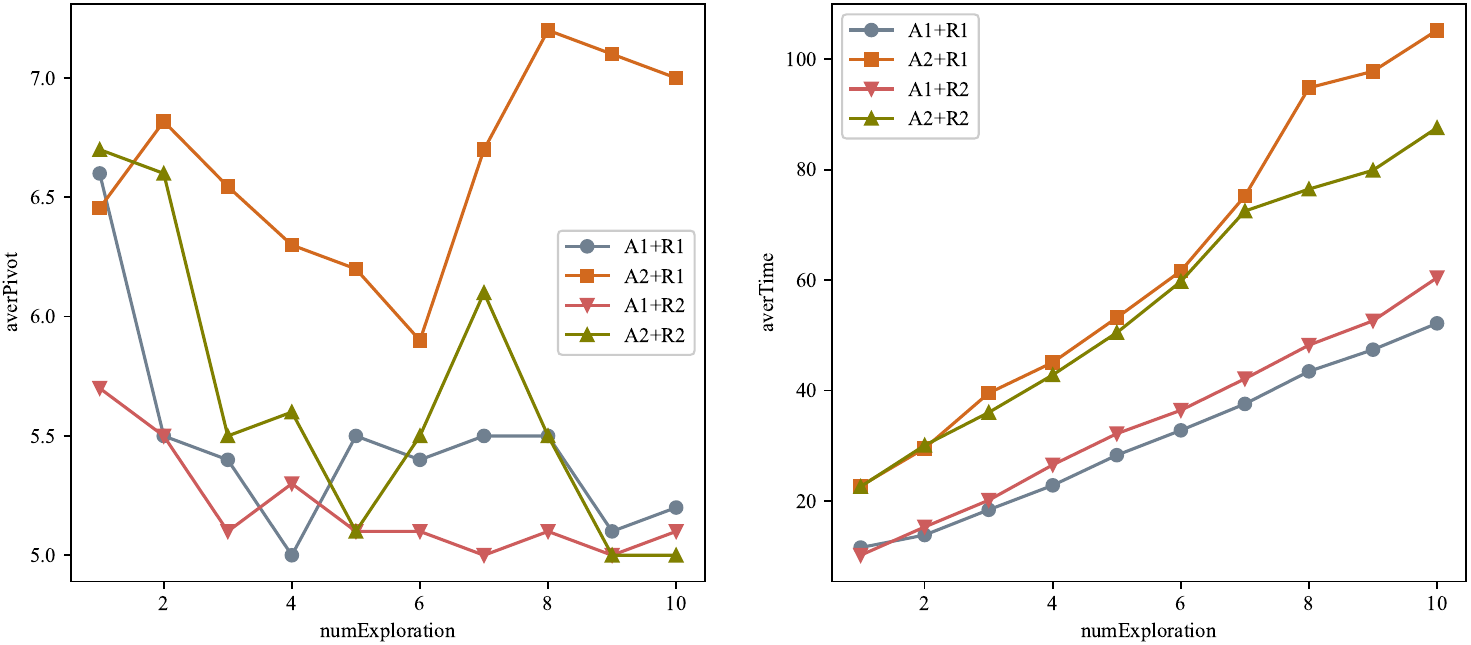}}
    \vspace{0.5em}
 
    \centering
    \subfigure[\label{fig:b}rand 300 $\times$ 300]
    {\includegraphics[width=1.00\columnwidth]{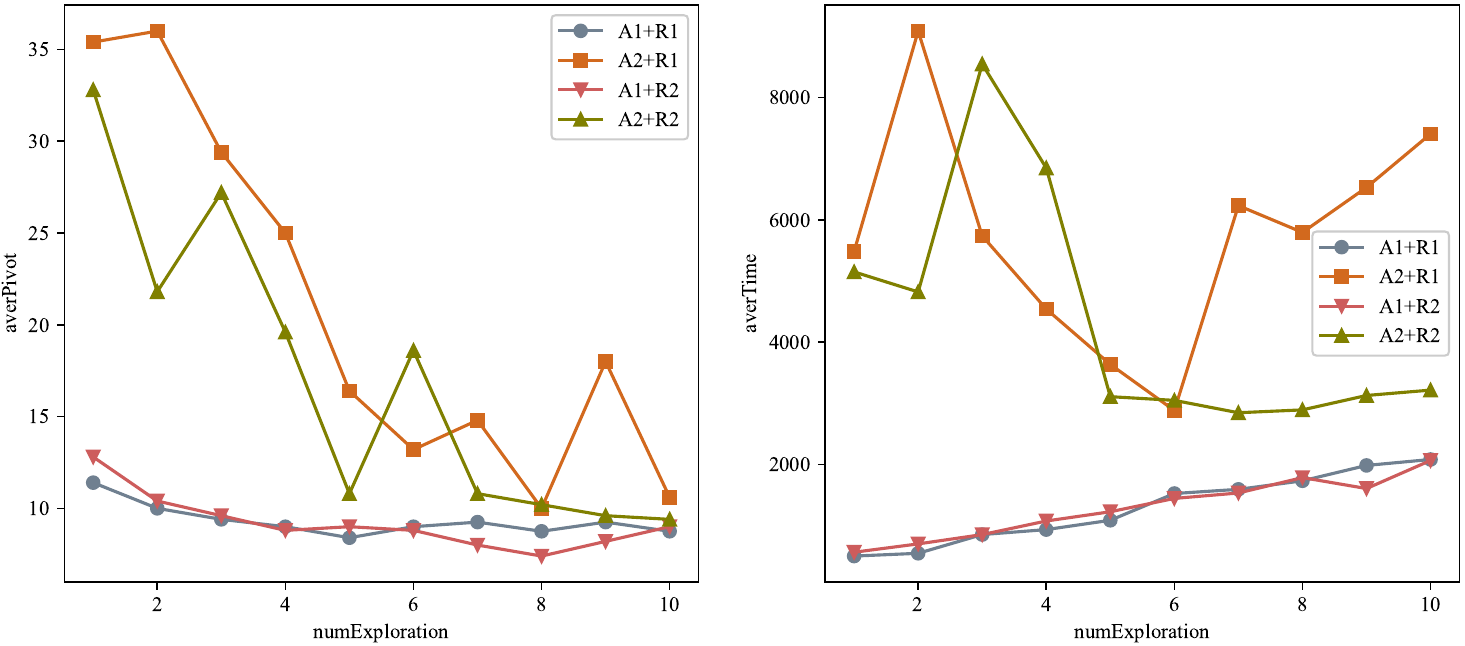}}
    \vspace{0.5em}
 
    \centering
    \subfigure[\label{fig:c}rand 232 $\times$ 504]
    {\includegraphics[width=1.00\columnwidth]{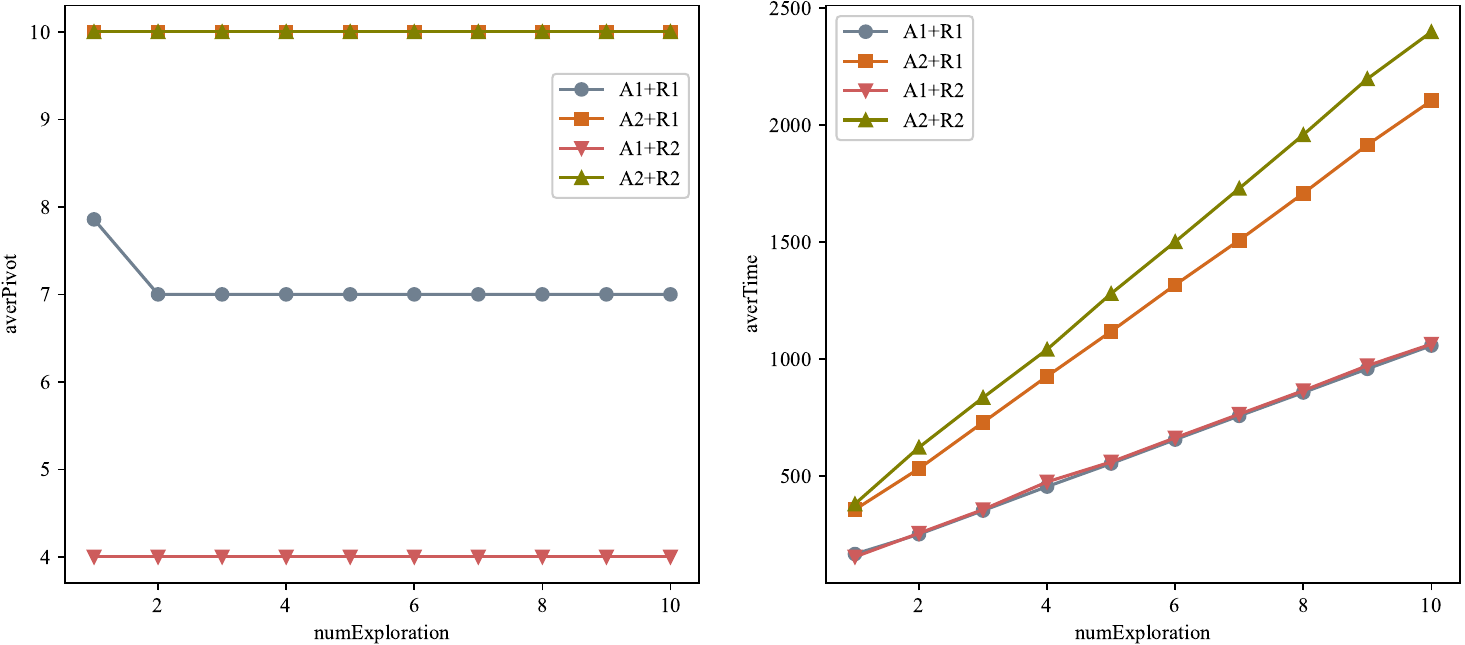}}
    \vspace{0.5em}
\end{figure}

\begin{figure}[htpb]
    \centering
    \subfigure[\label{fig:d}rand 332 $\times$ 187]
    {\includegraphics[width=1.00\columnwidth]{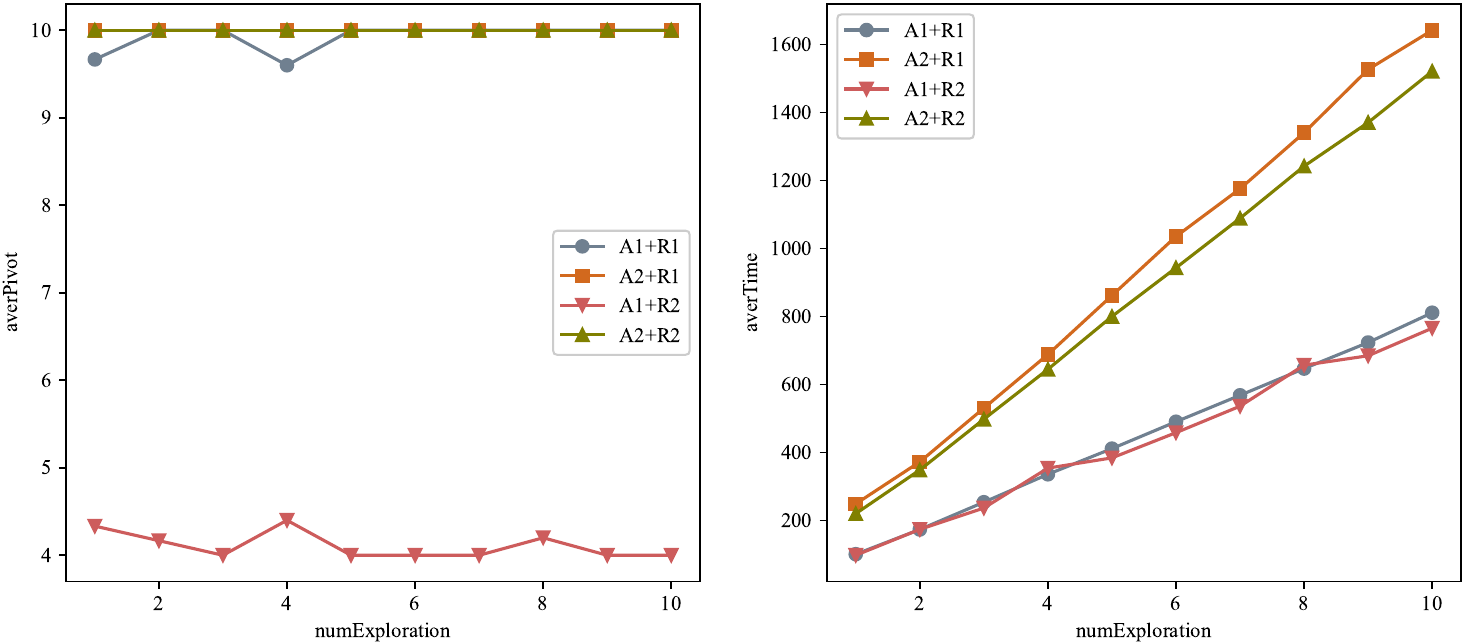}}
    \vspace{0.5em}
 
    \centering
    \subfigure[\label{fig:e}SC50A]
    {\includegraphics[width=1.00\columnwidth]{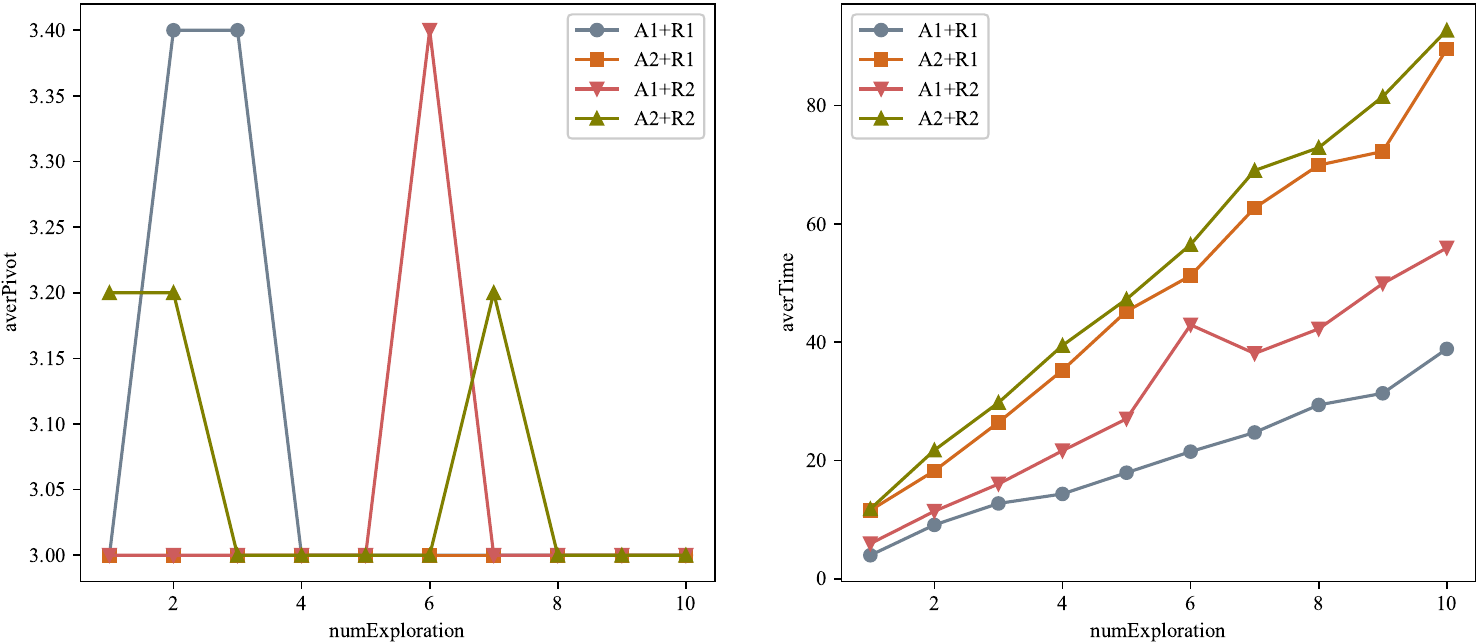}}
    \vspace{0.5em}
    
    \caption{Model comparison figure on five representative instances: rand 50 $\times$ 50, rand 232 $\times$ 504, rand 332 $\times$ 187 and SC50A. The X-axis represents the explorations, which are multiples of columns of the constraint matrix $A$. The Y-axis of the left figure represents the average pivot iterations. And the Y-axis of the right figure represents the average solution time.  }
    \label{fig:Model_ititalInstances}
\end{figure}

\paragraph{Implementation Details} In the experiment comparing the performance of the four models, we set the explorations from one to ten times the columns of the constraint matrix. Subsequently, we set the explorations to one and six times of columns because six times can provide stable exploration results, with number one having the shortest exploration time. We then solve each problem five times and obtain the optimal pivot iterations. In addition, we set $C=1 / \sqrt{2}$ in Formula (\ref{eq:equation_exploration_UCB}) and dynamically adjust $\alpha$ according to the solution state for Formula (\ref{eq:equation_exploration_Esoftmax}). In particular, we set $\alpha=0.3$ for $N_{explore}<{0.1 \times ColNum}$ and $\alpha=1$ for other cases, where $ColNum$ represents the number of columns of the constraint matrix. In addition, we use the first stage of the two-stage algorithm of the simplex method to find the initial feasible basis for the NETLIB instances.

\subsection{Estimation of Four RL Models}
In this section, we compare the quality of the four RL models proposed in Section 4.2 for the simplex method. We conduct thorough experiments on five representative instances, as illustrated in Figure~\ref{fig:Model_ititalInstances}. Each point represents the average pivot iterations obtained by solving the problem five times under the current conditions. This step reduces the influence of randomness. In addition, we set the upper limit of pivot iterations as $10$ to reduce unnecessary time without affecting the experimental results.

The left part of Figure~\ref{fig:Model_ititalInstances} shows the relationship between pivot iterations and explorations. We found that the number of pivot iterations obtained by A1+R2 is the best compared with the other three models. The right part of Figure~\ref{fig:Model_ititalInstances} shows the variation in the average solution time with the explorations. It can be seen that with the increase in explorations, the average solution time of each model also increases, which is caused by the exploration cost of the larger number of explorations. We also conclude that the time performance of A1+R1 and A1+R2 is better than those of A2+R1 and A2+R2. Consequently, combined with the performance of the number of pivot iterations, A1+R2 has certain advantages over all other models.

\subsection{Comparison with the Solver and Classical Pivot Rules}
In this section, we compare the MCTS rule with Python's linprog function and the other classic pivot rules. Table~\ref{tab:comparasion_rules} lists the results for the randomly generated square constraint matrices, Table~\ref{tab:comparasion_rules_general} presents the results for the general constraint matrices, and Table~\ref{tab:comparasion_rules_NETLIB} shows the experimental results for the NETLIB.

\begin{table*}
\centering
\caption{Comparison between classical pivot rules and the MCTS rule on random square constraint matrices. Red indicates the best results, and blue indicates suboptimal results.  }
\tiny
\resizebox{1.00\columnwidth}{!}{
\begin{tabular}{cc|ccccccc|c}
    \hline
    \textbf{row} & \textbf{col} & \textbf{Linprog} & \textbf{Danzig~\cite{DantzigRule}} & \textbf{Bland~\cite{BlandRule}} & \textbf{Steepest~\cite{steepestEdgeRule_1} } & \textbf{Greatest~\cite{GreatestImprovementRule}} & \textbf{Devex~\cite{DevexRule_2} } & \textbf{MCTS} & \textbf{MTime(s)}\\
    \hline
    \hline
    \multirow{9}{*}{300} & \multirow{9}{*}{300} & 341 & 18 & 63 & 14 & \textcolor{blue}{10} & 21 & \textcolor{red}{\textbf{7}} & 425 \\
    \cline{3-10}
    &&361 & 14 & 56 & 14 & \textcolor{blue}{13} & 22 & \textcolor{red}{\textbf{7}} & 409 \\
    \cline{3-10}
    && 331 & 14 & 49 & 12 & \textcolor{blue}{8} & 14 & \textcolor{red}{\textbf{6}} & 274 \\
    \cline{3-10}
    && 373 & 12 & 56 & 10 & \textcolor{blue}{8} & 16 & \textcolor{red}{\textbf{6}} & 287 \\
    \cline{3-10}
    && 318 & 16 & 28 & 13 & \textcolor{blue}{7} & 15 & \textcolor{red}{\textbf{6}} & 258 \\
    \cline{3-10}
    && 321 & 10 & 14 & \textcolor{blue}{8} & \textcolor{red}{\textbf{4}} & 14 & \textcolor{red}{\textbf{4}} & 128 \\
    \cline{3-10}
    && 324 & 17 & 67 & 14 & \textcolor{blue}{7} & 17 & \textcolor{red}{\textbf{5}} & 211 \\
    \cline{3-10}
    && 324 & 10 & 40 & 7 & \textcolor{blue}{6} & 14 & \textcolor{red}{\textbf{5}} & 130 \\
    \cline{3-10}
    && 329 & \textcolor{blue}{10} & 55 & 13 & \textcolor{red}{\textbf{5}} & 18 & \textcolor{red}{\textbf{5}} & 238 \\
    \hline
    \hline
    \multirow{9}{*}{400} & \multirow{9}{*}{400} & 492 & 21 & 37 & \textcolor{blue}{12} & \textcolor{blue}{12} & 29 & \textcolor{red}{\textbf{7}} & 356 \\
    \cline{3-10}
    && 475 & 22 & 113 & \textcolor{blue}{12} & \textcolor{red}{\textbf{9}} & 28 & \textcolor{red}{\textbf{9}} & 685 \\
    \cline{3-10}
    && 541 & 21 & 116 & 19 & \textcolor{blue}{15} & 32 & \textcolor{red}{\textbf{11}} & 793 \\
    \cline{3-10}
    && 449 & 19 & 73 & 15 & \textcolor{blue}{8} & 23 & \textcolor{red}{\textbf{6}} & 355 \\
    \cline{3-10}
    && 450 & 19 & 74 & 15 & \textcolor{blue}{11} & 30 & \textcolor{red}{\textbf{7}} & 510 \\
    \cline{3-10}
    && 473 & 16 & 70 & 13 & \textcolor{blue}{11} & 15 & \textcolor{red}{\textbf{6}} & 317 \\
    \cline{3-10}
    && 438 & 15 & 55 & \textcolor{blue}{9} & \textcolor{red}{\textbf{6}} & 24 & \textcolor{red}{\textbf{6}} & 275 \\
    \cline{3-10}
    && 480 & 15 & 73 & 12 & \textcolor{blue}{9} & 24 & \textcolor{red}{\textbf{8}} & 686 \\
    \cline{3-10}
    && 542 & 25 & 87 & 17 & \textcolor{blue}{10} & 32 & \textcolor{red}{\textbf{8}} & 553 \\
    \hline
    \hline
    \multirow{9}{*}{500} & \multirow{9}{*}{500} & 612 & 21 & 79 & \textcolor{blue}{14} & \textcolor{red}{\textbf{10}} & 32 &  \textcolor{red}{\textbf{10}} & 833 \\
    \cline{3-10}
    && 739 & 30 & 205 & \textcolor{blue}{20} & \textcolor{blue}{20} & 38 & \textcolor{red}{\textbf{15}} & 2633 \\
    \cline{3-10}
    && 688 & 22 & 93 & 18 & \textcolor{blue}{12} & 52 & \textcolor{red}{\textbf{11}} & 1064 \\
    \cline{3-10}
    && 559 & 23 & 74 & 13 & \textcolor{blue}{12} & 24 & \textcolor{red}{\textbf{8}} & 678 \\
    \cline{3-10}
    && 606 & 26 & 117 & 13 & \textcolor{blue}{12} & 28 & \textcolor{red}{\textbf{8}} & 664 \\
    \cline{3-10}
    && 528 & 22 & 44 & 11 & \textcolor{blue}{10} & 24 & \textcolor{red}{\textbf{6}} & 506 \\
    \cline{3-10}
    && 610 & 22 & 121 & 15 & \textcolor{blue}{10} & 27 & \textcolor{red}{\textbf{9}} & 966 \\
    \cline{3-10}
    && 598 & 21 & 58 & \textcolor{blue}{18} & \textcolor{red}{\textbf{10}} & 25 & \textcolor{red}{\textbf{10}} & 1395 \\
    \cline{3-10}
    && 607 & 26 & 187 & 17 & \textcolor{blue}{14} & 36 & \textcolor{red}{\textbf{11}} & 1383 \\
    \hline
\end{tabular}}
\label{tab:comparasion_rules}
\end{table*}

\subsubsection{Comparison Results on Random instances}
The instances in Table~\ref{tab:comparasion_rules} are random square matrices. However, Table~\ref{tab:comparasion_rules_general} presents a comparison of the random rows and columns of the constraint matrices. Red indicates the best results, and blue indicates suboptimal results. The first two columns represent rows and columns of the constraint matrix. The third column to the ninth column respectively represent the minimum number of pivot iterations of Python's linprog function, the Danzig rule~\cite{DantzigRule}, the Bland rule~\cite{BlandRule}, the steepest-edge rule~\cite{steepestEdgeRule_1}, the greatest improvement rule~\cite{GreatestImprovementRule}, the devex rule~\cite{DevexRule_2} and our MCTS rule. The initial feasible basis of Python's linprog function is determined by its own setting, while others are based on the method proposed in Section 5.1. The last column is execution time of the MCTS rule.

Tables~\ref{tab:comparasion_rules} and \ref{tab:comparasion_rules_general} show that the number of pivot iterations obtained by the MCTS rule is superior to that obtained by the classical pivot rules in all general instances. In terms of square instances, the pivot iterations found by the MCTS rule were only 54.55\% of the minimum iterations of other popular pivot rules. In addition, for other random dimension instances, our result was only 55.56\% of the others’ best results.

We conclude that the results of the MCTS rule are not limited to input instances. It performed best on all randomly generated problems because the MCTS rule selects the entering basis variable by exploring and evaluating the entire feasible action space rather than providing a fixed rule based on certain specific features. However, the number of pivot iterations obtained by other classical methods cannot exceed the result of the MCTS rule. We first proposed an efficient and generalized method for determining the minimum number of pivot iterations of the simplex method. Furthermore, for the first time, this method provides the best label design for pivot rules based on supervised learning.

{
\tiny
\setlength{\tabcolsep}{.5mm}
\setlength{\LTcapwidth}{4.7in}
\begin{longtable}{cc|ccccccc|c}
\caption{Comparison between classical pivot rules and the MCTS rule on the general constraint matrices. Red indicates the best results, and blue indicates suboptimal results. }
    \\
    \hline
    \textbf{row} & \textbf{col} & \textbf{Linprog} & \textbf{Danzig~\cite{DantzigRule}} & \textbf{Bland~\cite{BlandRule}} & \textbf{Steepest~\cite{steepestEdgeRule_1} } & \textbf{Greatest~\cite{GreatestImprovementRule}} & \textbf{Devex~\cite{DevexRule_2} } & \textbf{MCTS} & \textbf{MTime(s)} \\
    \hline
    117 & 123 & 187 & 8 & 30 & 7 & \textcolor{blue}{5} & 8 & \textcolor{red}{\textbf{4}} & 35 \\
    \hline
    139 & 199 & 191 & 9 & 24 & \textcolor{blue}{7} & \textcolor{red}{\textbf{5}} & 10 & \textcolor{red}{\textbf{5}} & 59 \\
    \hline
    399 & 324 & 571 & 25 & 1000+ & \textcolor{blue}{14} & 17 & 44 & \textcolor{red}{\textbf{9}} & 1107 \\
    \hline
    243 & 245 & 267 & 10 & 30 & \textcolor{blue}{8} & \textcolor{red}{\textbf{5}} & 16 & \textcolor{red}{\textbf{5}} & 168 \\
    \hline
    471 & 311 & 529 & 8 & 46 & 7 & \textcolor{blue}{5} & 9 & \textcolor{red}{\textbf{4}} & 191 \\
    \hline
    209 & 318 & 324 & \textcolor{blue}{9} & 64 & \textcolor{blue}{9} & \textcolor{blue}{9} & \textcolor{blue}{9} & \textcolor{red}{\textbf{7}} & 288 \\
    \hline
    374 & 293 & 478 & 11 & 22 & 8 & \textcolor{blue}{7} & 16 & \textcolor{red}{\textbf{6}} & 351 \\
    \hline
    470 & 766 & 765 & 485 & 64 & 10 & \textcolor{blue}{8} & 21 & \textcolor{red}{\textbf{6}} & 1133 \\
    \hline
    49 & 647 & 105 & 1000+ & 1000+ & \textcolor{blue}{7} & \textcolor{red}{\textbf{4}} & 21 & \textcolor{red}{\textbf{4}} & 427 \\
    \hline
    207 & 327 & 259 & \textcolor{blue}{6} & 34 & \textcolor{blue}{6} & \textcolor{red}{\textbf{4}} & 11 & \textcolor{red}{\textbf{4}} & 145 \\
    \hline
    151 & 419 & 260 & 1000+ & 36 & 13 & \textcolor{blue}{11} & 22 & \textcolor{red}{\textbf{7}} & 735 \\
    \hline
    782 & 436 & 871 & \textcolor{blue}{11} & 43 & 12 & \textcolor{red}{\textbf{5}} & 25 & \textcolor{red}{\textbf{5}} & 917 \\
    \hline
    363 & 699 & 418 & 8 & 41 & 7 & \textcolor{blue}{5} & 23 & \textcolor{red}{\textbf{4}} & 555 \\
    \hline
    587 & 382 & 635 & \textcolor{blue}{8} & 33 & 12 & \textcolor{red}{\textbf{4}} & 9 & \textcolor{red}{\textbf{4}} & 416 \\
    \hline
    565 & 761 & 594 & 6 & 10 & \textcolor{blue}{5} & \textcolor{red}{\textbf{2}} & 8 & \textcolor{red}{\textbf{2}} & 214 \\
    \hline
    232 & 504 & 263 & \textcolor{blue}{6} & 110 & \textcolor{blue}{6} & \textcolor{red}{\textbf{4}} & 7 & \textcolor{red}{\textbf{4}} & 152 \\
    \hline
    215 & 202 & 358 & 18 & 58 & 13 & \textcolor{blue}{8} & 22 & \textcolor{red}{\textbf{7}} & 250 \\
    \hline
    278 & 525 & 373 & 8 & 52 & \textcolor{blue}{7} & 10 & 14 &  \textcolor{red}{\textbf{5}} & 648 \\
    \hline
    143 & 58 & 184 & 6 & 15 & 6 & \textcolor{blue}{4} & 6 &  \textcolor{red}{\textbf{3}} & 26 \\
    \hline
    323 & 286 & 378 & \textcolor{blue}{5} & 17 & \textcolor{blue}{5} & \textcolor{red}{\textbf{4}} & 8 &  \textcolor{red}{\textbf{4}} & 191 \\
    \hline
    133 & 242 & 192 & 18 & 1000+ & 9 & \textcolor{blue}{7} & 21 &  \textcolor{red}{\textbf{5}} & 282 \\
    \hline
    96 & 482 & 179 & 9 & 46 & \textcolor{blue}{8} & \textcolor{red}{\textbf{5}} & 11 &  \textcolor{red}{\textbf{5}} & 397 \\
    \hline
    203 & 200 & 223 & 8 & 12 & \textcolor{blue}{6} & \textcolor{red}{\textbf{3}} & 7 &  \textcolor{red}{\textbf{3}} & 84 \\
    \hline
    678 & 154 & 757 & 10 & 14 & \textcolor{blue}{7} & \textcolor{red}{\textbf{3}} & 8 &  \textcolor{red}{\textbf{3}} & 267 \\
    \hline
    739 & 142 & 866 & \textcolor{blue}{7} & 18 & \textcolor{blue}{7} & 8 & 11 &  \textcolor{red}{\textbf{5}} & 520 \\
    \hline
    739 & 730 & 841 & 11 & 33 & 10 & \textcolor{red}{\textbf{4}} & \textcolor{blue}{7} &  \textcolor{red}{\textbf{4}} & 888 \\
    \hline
    240 & 234 & 268 & 7 & 20 & 7 & \textcolor{blue}{5} & 10 &  \textcolor{red}{\textbf{4}} & 103 \\
    \hline
    767 & 158 & 885 & 9 & 1000+ & \textcolor{blue}{7} & \textcolor{blue}{7} & 19 &  \textcolor{red}{\textbf{5}} & 623 \\
    \hline
    730 & 467 & 773 & 6 & 21 & 7 & \textcolor{blue}{5} & 6 &  \textcolor{red}{\textbf{3}} & 584 \\
    \hline
    434 & 258 & 498 & 7 & 30 & 7 & \textcolor{blue}{5} & 23 &  \textcolor{red}{\textbf{4}} & 231 \\
    \hline
    625 & 521 & 739 & 9 & 26 & 10 & \textcolor{blue}{9} & 19 &  \textcolor{red}{\textbf{5}} & 732 \\
    \hline
    332 & 187 & 352 & \textcolor{blue}{9} & 24 & \textcolor{blue}{9} & \textcolor{red}{\textbf{4}} & 13 &  \textcolor{red}{\textbf{4}} & 102 \\
    \hline
    561 & 628 & 585 & \textcolor{blue}{7} & 13 & 9 & \textcolor{red}{\textbf{4}} & 16 &  \textcolor{red}{\textbf{4}} & 471 \\
    \hline
    587 & 774 & 610 & 5 & 17 & 5 & \textcolor{blue}{4} & 13 &  \textcolor{red}{\textbf{3}} & 302 \\
    \hline
    84 & 108 & 113 & 15 & 12 & 7 & \textcolor{blue}{5} & 14 &  \textcolor{red}{\textbf{4}} & 35 \\
    \hline
\label{tab:comparasion_rules_general}\\
\end{longtable}
}

\subsubsection{Comparison Results on NETLIB}
Table~\ref{tab:comparasion_rules_NETLIB} presents the comparison results of the MCTS rule with other classical pivot rules on NETLIB. Red indicates the best results and blue indicates suboptimal results. The first column represents the name of instances. The second to seventh columns are respectively the minimum pivot iterations of Python’s linprog function, the Danzig rule~\cite{DantzigRule}, the Bland rule~\cite{BlandRule}, the steepest-edge rule~\cite{steepestEdgeRule_1}, the greatest improvement rule~\cite{GreatestImprovementRule} and our MCTS rule. The initial feasible basis of Python's linprog function is determined by its own setting, while others are based on the method proposed in Section 5.1. The last column is execution time of the MCTS rule.

It is easy to conclude that the MCTS rule yields the least number of pivot iterations far less than others on all the instances listed, especially for the problem ADLITTLE. The greatest improvement rule has the least number of pivot iterations among the classical rules. In contrast, our method achieves only 49.06\% of its pivot iterations. Moreover, compared with the Danzig rule, the MCTS rule gets less than 2.6\% of its pivot iterations.

\begin{table*}[thpb]
\centering
\caption{Comparison between classical pivot rules and the MCTS rule on NETLIB. Red indicates the best results, and blue indicates suboptimal results. }
\tiny
\setlength{\tabcolsep}{.5mm}
\resizebox{1.00\columnwidth}{!}{
\begin{tabular}{c|cccccc|c}
    \hline
    \textbf{Problem} & \textbf{Linprog} & \textbf{Danzig~\cite{DantzigRule}} & \textbf{Bland~\cite{BlandRule}} & \textbf{Steepest~\cite{steepestEdgeRule_1} } & \textbf{Greatest~\cite{GreatestImprovementRule}} & \textbf{MCTS} & \textbf{MTime(s)}\\
    \hline
    AFIRO & 27 & \textcolor{red}{\textbf{0}} & \textcolor{red}{\textbf{0}} & \textcolor{red}{\textbf{0}} & \textcolor{red}{\textbf{0}} & \textcolor{red}{\textbf{0}} & 0 \\
    \hline
    ADLITTLE & 152 & 1000+ & 223 & 60 & \textcolor{blue}{53} & \textcolor{red}{\textbf{26}} & 836 \\
    \hline
    BLEND & 178 & 29 & 88 & 35 & \textcolor{blue}{21} & \textcolor{red}{\textbf{15}} & 118 \\
    \hline
    SC50A & 56 & 5 & 8 & 5 & \textcolor{blue}{4} & \textcolor{red}{\textbf{3}} & 4 \\
    \hline
    SC50B & 59 & \textcolor{red}{\textbf{6}} & 9 & \textcolor{blue}{7} & \textcolor{red}{\textbf{6}} & \textcolor{red}{\textbf{6}} & 7 \\
    \hline
    SC105 & 135 & 14 & 26 & \textcolor{blue}{10} & 16 & \textcolor{red}{\textbf{7}} & 69 \\
    \hline
    SCAGR7 & 228 & 57 & 101 & \textcolor{blue}{43} & 54 & \textcolor{red}{\textbf{40}} & 1488 \\
    \hline
    SHARE2B & 260 & 47 & 172 & 38 & \textcolor{blue}{29} & \textcolor{red}{\textbf{19}} & 612 \\
    \hline
\end{tabular}}
\label{tab:comparasion_rules_NETLIB}
\end{table*}

Although the solution time of the MCTS rule is longer than other algorithms, this is still consistent with our contribution. We aim to determine the optimal pivot iterations and all corresponding pivot paths for the input instance. Furthermore, we provide the best supervision labels for the simplex method. Additionally, in Section 7, we present two methods from the perspective of CPU and GPU to improve the efficiency of collecting supervision labels. Thus, this method is more applicable to super large-scale problems.

\subsection{Comparisons with all Current Pivot Rules Based on ML}
In this section, we compare our method with two recently proposed machine-learning methods. It is found that the minimum pivot iteration MCTS rule finds is much better than the other two methods.

In the first method~\cite{1}, the supervised learning method DeepSimplex~\cite{1} performs better than the unsupervised learning method. Therefore, we only compare DeepSimplex~\cite{1}, which gets the best performance of the Dantzig rule and steepest-edge rule as the supervised signal. Thus, the results do not exceed those of these two methods. As Tables ~\ref{tab:comparasion_rules} and \ref{tab:comparasion_rules_general} show, the worst number of pivot iterations of the MCTS rule in all instances is $80\%$ of Danzig's. And in the best case, the number of pivot iterations is less than $0.4\%$ of Danzig's. Compared with the steepest-edge rule, the worst pivot iterations of the MCTS rule in all instances is only $77.78\%$. Moreover, in the best case, the pivot iteration is only $27.27\%$ of the steepest-edge rule. Additionally, Table~\ref{tab:comparasion_rules_NETLIB} shows that the MCTS rule is significantly better than the Dantzig and steepest-edge rules on NETLIB. Therefore, we conclude that the MCTS rule can obtain better pivot iterations than DeepSimplex~\cite{1}. 

The second method~\cite{2} is learned in a supervised manner with the best label of Dantzig's rule, the hybrid (DOcplex's default) rule, the greatest improvement rule, the steepest-edge rule, and the devex rule. Especially remarkable is that the MCTS rule provided more effective labels with minimum pivot iterations. From the experimental results of their article~\cite{2}, we know that the best result of the average number of pivot iterations is $99.54\%$ of the number of steepest-edge rule. In contrast, the performance of the MCTS rule on the worst instance is $77.78\%$ of the steepest-edge rule, as shown in Tables~\ref{tab:comparasion_rules} and \ref{tab:comparasion_rules_general}. Moreover, the best result obtained by the MCTS rule is only $33.33\%$ of the pivot iterations of the steepest-edge rule. Furthermore, in NETLIB, the pivot iterations yielded by the MCTS rule can even reach 43.33\% of the steepest-edge rule, which is far less than 99.54\%.

\begin{table}[t]
\centering
\caption{Pivoting paths with the minimum number of pivot iterations on five representative instances. The pivoting path is an ordered sequence of entering basis variables.}
\tiny
\setlength{\belowcaptionskip}{15pt}
\resizebox{1.00\columnwidth}{!}{
\begin{tabular}{c|cccc}
    \hline
    \textbf{Problem} & \textbf{Index} & \textbf{Pivot paths} & \textbf{Pivot iterations} & \textbf{Objective value} \\
    \hline
    \hline
    \multirow{6}{*}{SC50B} & 1 & $\left[1, 24, 36, 12, 23, 35 \right]$ & 6 & 70.000 \\
    \cline{2-5}
    & 2 & $\left[1, 35, 12, 24, 36, 23 \right]$ & 6 & 70.000 \\
    \cline{2-5}
    & 3 & $\left[1, 35, 23, 36, 12, 24 \right]$ & 6 & 70.000 \\
    \cline{2-5}
    & 4 & $\left[1, 35, 24, 12, 36, 23 \right]$ & 6 & 70.000 \\
    \cline{2-5}
    & 5 & $\left[12, 35, 1, 24, 36, 23 \right]$ & 6 & 70.000 \\
    \cline{2-5}
    & 6 & $\left[23, 35, 1, 36, 12, 24 \right]$ & 6 & 70.000 \\
    \hline
    \hline
    \multirow{2}{*}{rand 232*504} & 1 & $\left[52, 55, 358, 220 \right]$ & 4 & 22.287 \\
    \cline{2-5}
    & 2 & $\left[55, 52, 358, 220 \right]$ & 4 & 22.287 \\
    \hline
    \hline
    \multirow{5}{*}{rand 434*258} & 1 & $\left[14, 155, 223, 150 \right]$ & 4 & 5.509 \\
    \cline{2-5}
    & 2 & $\left[14, 155, 150, 223 \right]$ & 4 & 5.509 \\
    \cline{2-5}
    & 3 & $\left[150, 155, 223, 14 \right]$ & 4 & 5.509 \\
    \cline{2-5}
    & 4 & $\left[155, 150, 223, 14 \right]$ & 4 & 5.509 \\
    \cline{2-5}
    & 5 & $\left[155, 223, 150, 14 \right]$ & 4 & 5.509 \\
    \hline
    \hline
    \multirow{3}{*}{rand 50*50} & 1 & $\left[ 1,48,28,19,21 \right]$ & 5 & 3.482 \\
    \cline{2-5}
    & 2 & $\left[ 28,48,1,19,21 \right]$ & 5 & 3.482 \\
    \cline{2-5}
    & 3 & $\left[ 48,28,1,19,21 \right]$ & 5 & 3.482 \\
    \hline
    \hline
    \multirow{6}{*}{rand 300*300} & 1 & $\left[26, 294, 42, 143, 56, 116, 263 \right]$ & 7 & 5.915 \\
    \cline{2-5}
    & 2 & $\left[26, 294, 42, 116, 56, 143, 263 \right]$ & 7 & 5.915 \\
    \cline{2-5}
    & 3 & $\left[56, 26, 42, 294, 263, 116, 143 \right]$ & 7 & 5.915 \\
    \cline{2-5}
    & 4 & $\left[143, 26, 42, 294, 116, 56, 263 \right]$ & 7 & 5.915 \\
    \cline{2-5}
    & 5 & $\left[143, 26, 56, 294, 42, 116, 263 \right]$ & 7 & 5.915 \\
    \cline{2-5}
    & 6 & $\left[294, 26, 42, 143, 56, 116, 263 \right]$ & 7 & 5.915 \\
    \hline
\end{tabular}}
\label{tab:sequences}
\end{table}

\begin{figure*}
\centering
\includegraphics[width=.7\columnwidth]{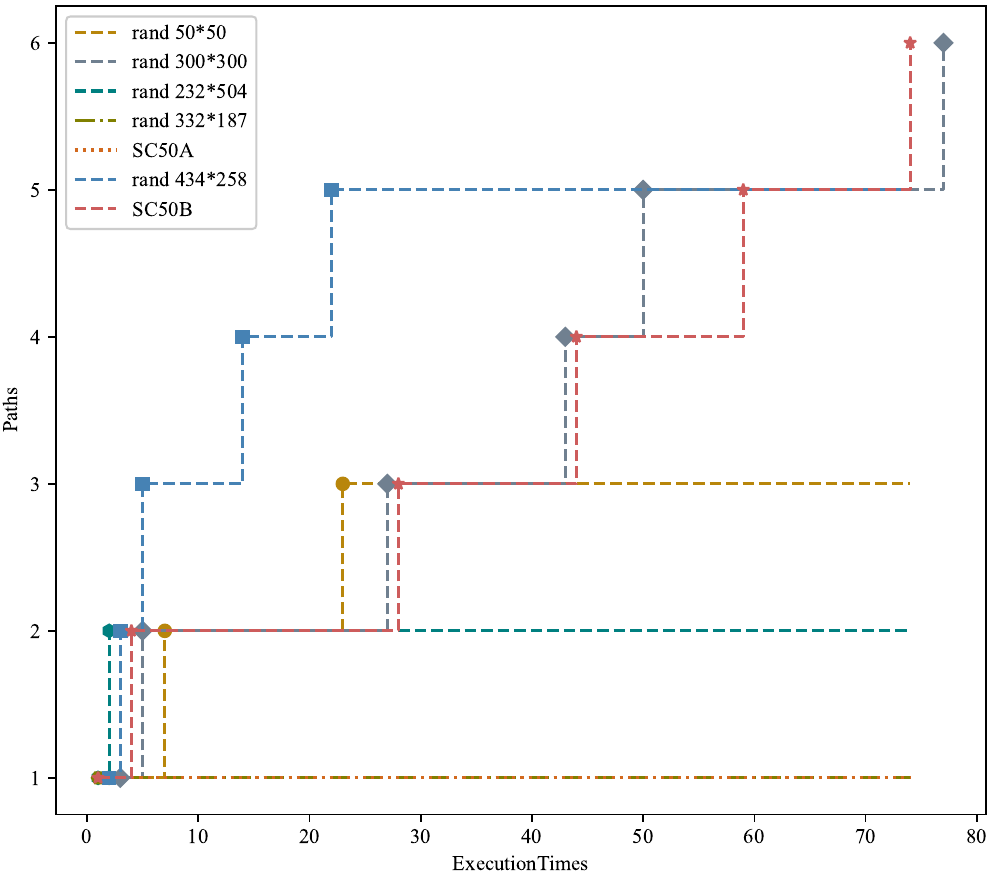} 
\caption{Multiple paths found vary with the number of algorithm executions. The X-axis represents the number of algorithm executions, and the Y-axis represents the different pivot paths currently found. }
\label{fig:multiple_paths}
\end{figure*}

\subsection{Findings of Multiple Pivot Paths}
Providing multiple pivot paths for the simplex method is the result of taking advantage of the randomness of the MCTS rule. Table~\ref{tab:sequences} shows multiple optimal pivoting paths found for five representative problems, which cannot be yielded by previous methods. Different pivoting paths are the optimal pivoting sequences with minimum pivot iterations. Additionally, Figure~\ref{fig:multiple_paths} shows the relationship between the number of different pivot paths found and the algorithm executions on several representative instances mentioned above. We use highlighted points to mark the executions of a newly found pivot path. It can be concluded that under the initial executions, the proposed MCTS rule can find some paths.  Furthermore, with the increase of algorithm executions, the number of found paths also increases.

\begin{figure*}[htbp]
    \centering
    \subfigure[\label{fig:a}GlobalFigure]  
    {  
        \begin{minipage}[t]{0.50\linewidth}  
        \centering
        \includegraphics[width=1.00\columnwidth]{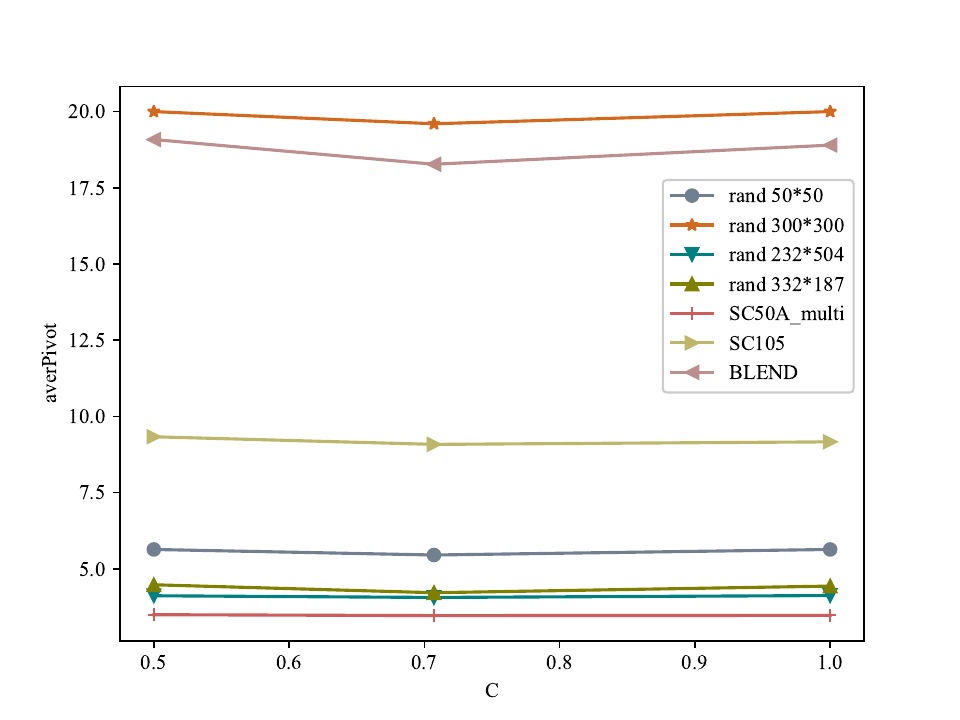}
        \end{minipage}%
    }%
    \subfigure[\label{fig:b}LocalFigure]  
    {
        \begin{minipage}[t]{0.5\linewidth}  
        \centering
        \includegraphics[width=1.00\columnwidth]{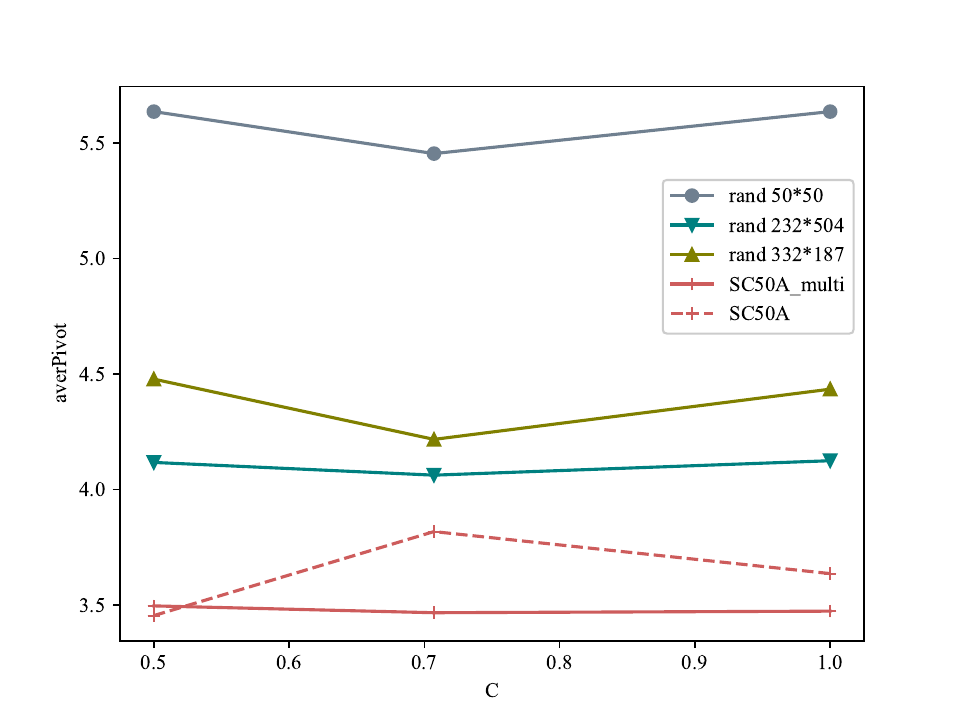}
        \end{minipage}%
    }%
\caption{Relationship between average pivot iterations and parameter $C$. The X-axis represents the value of parameter $C$, and the Y-axis represents the average pivot iterations found. On the left is the overall effect of several representative instances. On the right is an enlarged representation of the bottom four lines in the left figure. The SC50A$\_$multi represents an increase in the executions of SC50A by five times than before to calculate the average value. }  
\label{fig:C}  
\end{figure*}

\begin{figure*}[htbp]
    \centering
    \subfigure[\label{fig:a}ExploreNum = 1 $\times$ ColNum]  
    {  
        \begin{minipage}[t]{0.50\linewidth}  
        \centering
        \includegraphics[width=1.00\columnwidth]{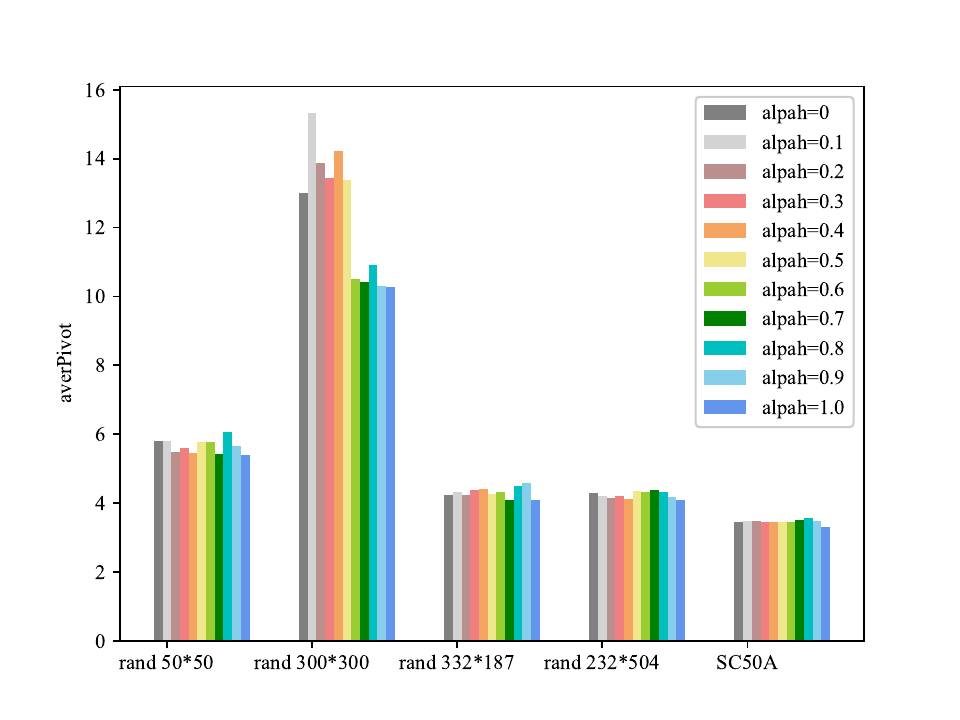}
        \end{minipage}%
    }%
    \subfigure[\label{fig:b}ExploreNum = 0.5 $\times$ ColNum]  
    {
        \begin{minipage}[t]{0.5\linewidth}  
        \centering
        \includegraphics[width=1.00\columnwidth]{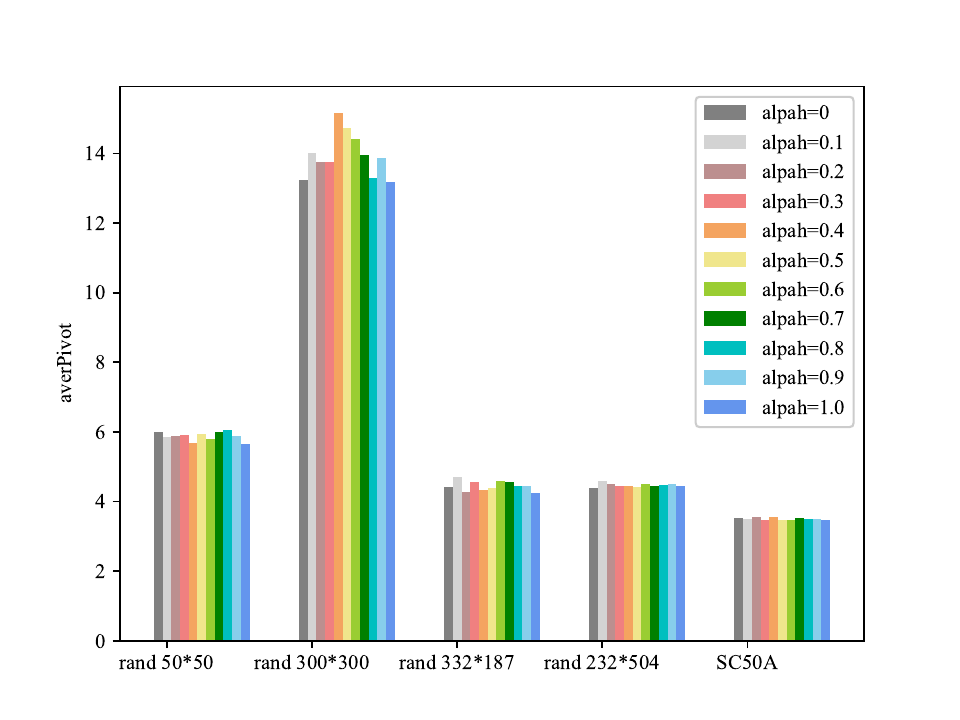}
        \end{minipage}%
    }%
    
    \centering
    \subfigure[\label{fig:c}ExploreNum = 0.4 $\times$ ColNum]  
    {  
        \begin{minipage}[t]{0.50\linewidth}  
        \centering
        \includegraphics[width=1.00\columnwidth]{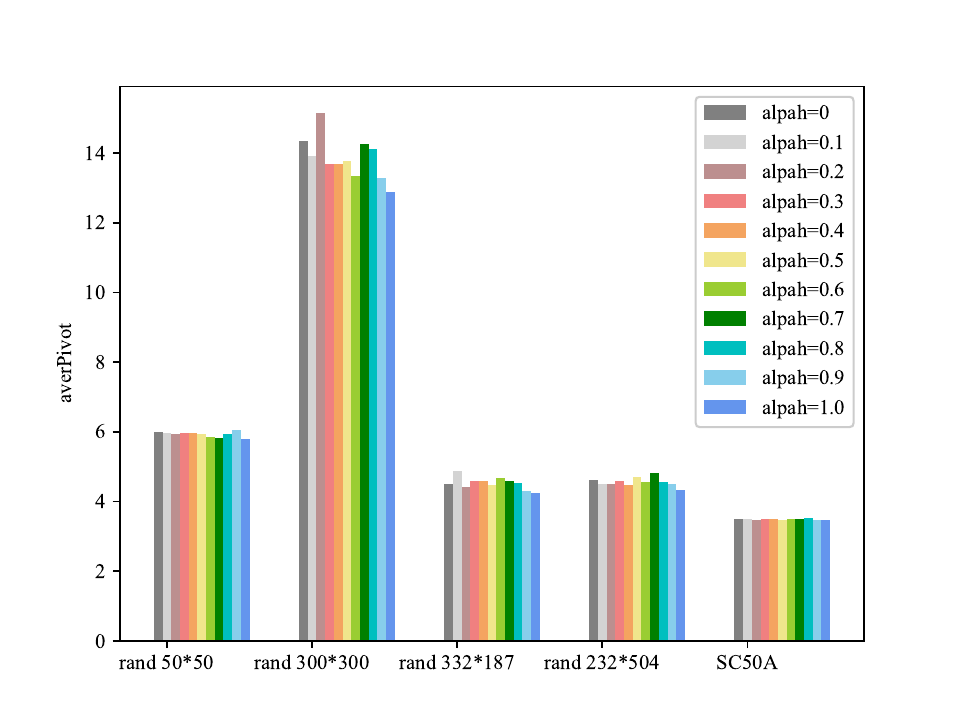}
        \end{minipage}%
    }%
    \subfigure[\label{fig:d}ExploreNum = 0.3 $\times$ ColNum]  
    {
        \begin{minipage}[t]{0.5\linewidth}  
        \centering
        \includegraphics[width=1.00\columnwidth]{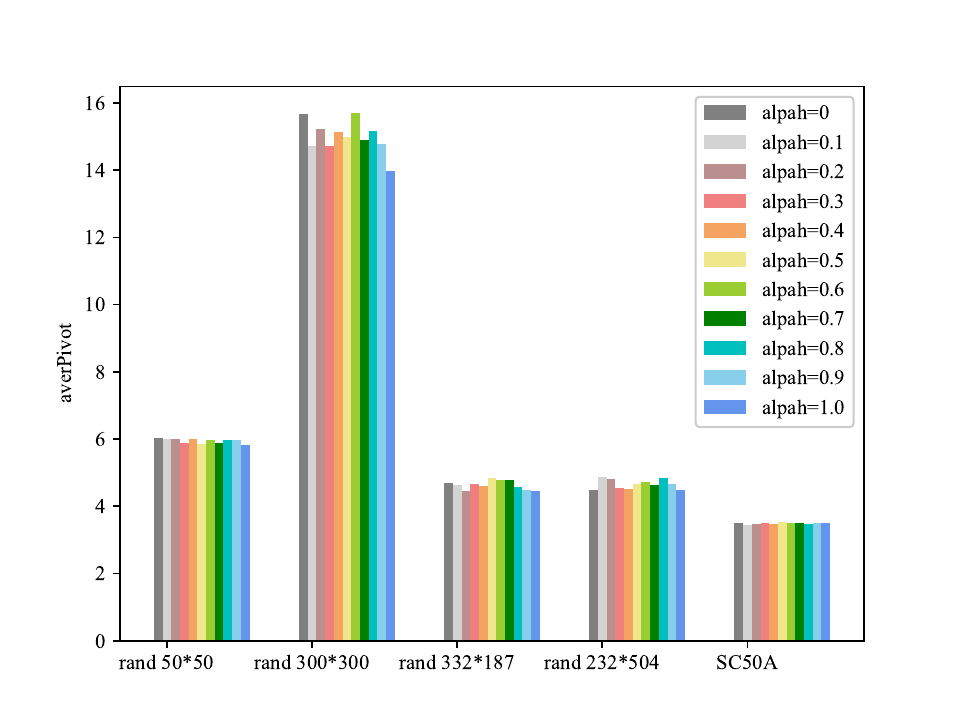}
        \end{minipage}%
    }
    
    \centering
    \subfigure[\label{fig:c}ExploreNum = 0.2 $\times$ ColNum]  
    {  
        \begin{minipage}[t]{0.50\linewidth}  
        \centering
        \includegraphics[width=1.00\columnwidth]{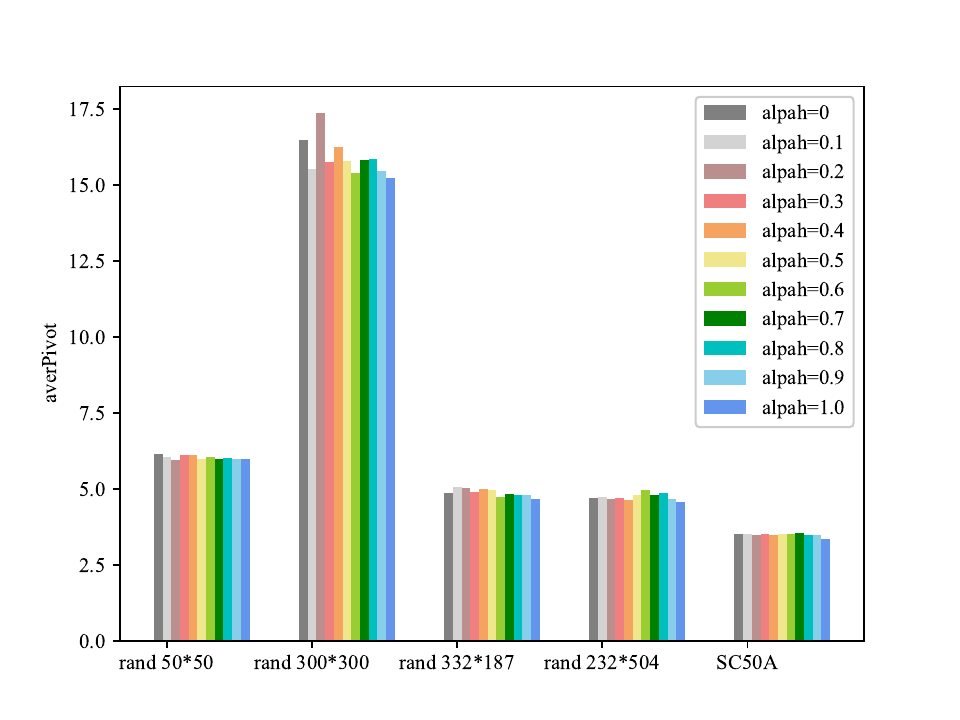}
        \end{minipage}%
    }%
    \subfigure[\label{fig:d}ExploreNum = 0.1 $\times$ ColNum]  
    {
        \begin{minipage}[t]{0.5\linewidth}  
        \centering
        \includegraphics[width=1.00\columnwidth]{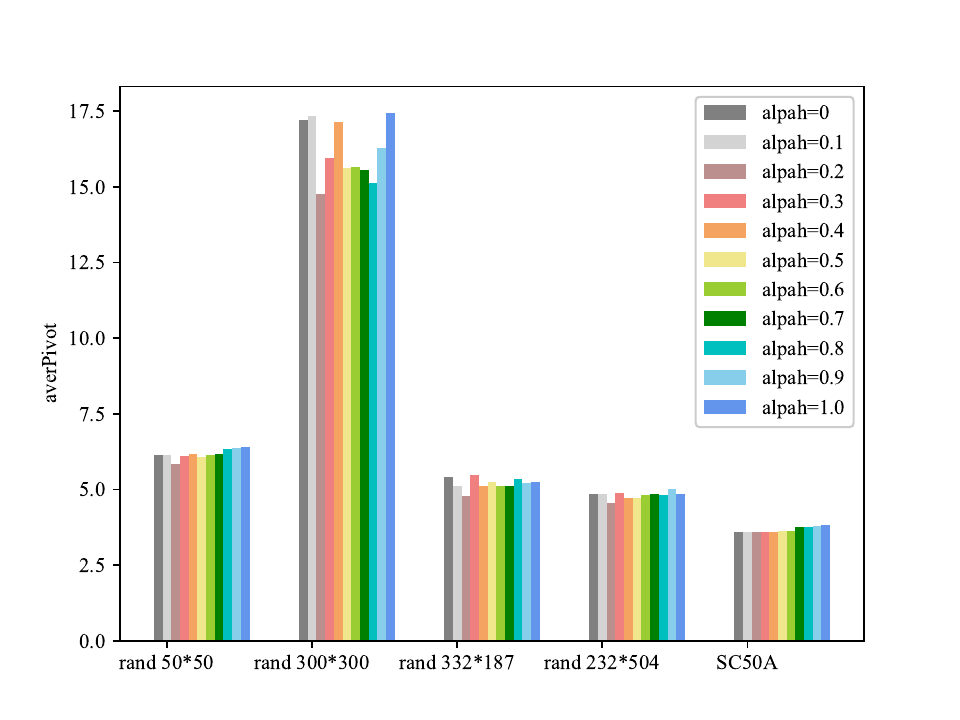}
        \end{minipage}%
    }%
\caption{Relationship between the average number of pivot iterations and $\alpha$ under different initial explorations. The X-axis represents different problems, and the Y-axis represents the average pivot iterations. }  
\label{fig:Alpha}  
\end{figure*}

\subsection{Ablation Study}
We compare the influence of different $C$ and $\alpha$ values on the average pivot iterations for several representative instances. Each point is the average result of executing the algorithm five times. As the GlobalFigure in Figure~\ref{fig:C} indicates, the empirical value of $C$ can lead to the least pivot iterations, except for SC50A. However, we conduct in-depth experiments and find that SC50A can provide fewer pivot iterations when the total number of executions increases. Therefore, we believe that the empirical value of $C$ is reasonable for the MCTS rule in terms of overall performance.

Figure~\ref{fig:Alpha} shows the relationship between $\alpha$ and the pivot iteations for different initial explorations. Formula (\ref{eq:equation_exploration_select}) aims to relax the $max$ operator owing to the imprecisely estimated value in the early stage of exploration. Therefore, we conduct sufficient experiments on explorations of the 1-, 0.5-, 0.4-, 0.3-, 0.2-, and 0.1-times columns of the constraint matrix of the instances to be solved. Formula (\ref{eq:equation_exploration_select}) is effective when explorations are less than or equal to $0.1$ times columns. Furthermore, when the number of explorations is $0.1$ times columns, $\alpha$ achieves a consistently good effect with a value of $0.3$. Therefore, we set dynamically adjusted $\alpha$ for the MCTS rule. When the number of explorations is less than or equal to $0.1$ times columns, $\alpha$ is set to $0.3$. For the other cases, $\alpha$ was set to $1$.

\section{Conclusion}
Based on the proposed SimplexPseudoTree structure and reinforcement learning model, the MCTS rule can determine all the shortest pivot paths of the simplex method. In addition, our method provides the best supervised label-setting method for the simplex method based on supervised learning. The MCTS rule can evaluate the pros and cons of entering basis variables individually, significantly reducing the exploration space for combinatorial optimization problems. Therefore, the proposed method can find the minimum pivot iterations and provide a method to find multiple shortest pivot paths. This idea can also be used to find multiple optimal solutions for other combinatorial optimization problems that can be modeled as imitative tree structures. Furthermore, we prove that the MCTS rule can find polynomial pivot iterations when the number of vertices in the feasible region is $C_n^m$. The complete theory and comprehensive experiments demonstrate that the MCTS rule can find multiple optimal pivot rules.

\section{Further Work}
The multiple pivot paths determined by the MCTS rule can be used to construct flexible labels for the simplex method. Therefore, we can design the supervised learning method of the optimal pivot rule for the simplex method of linear programming. Furthermore, deep learning can be used to construct more efficient and time-effective pivot rules. In this manner, we can improve the redundancy of the traditional pivot rule and the low time efficiency of the MCTS rule.

Additionly, we introduce two implementation techniques to improve the time efficiency of the proposed method. These techniques are introduced from the perspectives of the CPU and GPU. Both methods are designed to solve the time-consuming process of sequential execution in the exploration stage. First, rewriting CUDA allows several explorations to be performed simultaneously. Thus, the time efficiency is reduced by dozens or even hundreds of times. In addition to using GPU computing by rewriting CUDA, the implementation of multithreading provides a method to improve the time efficiency of CPU devices. $N_{explore}$ explorations can be divided into $N_{explore}/N_{threads}$ groups by grouping, where $N_{explore}$ represents the number of explorations and $N_{threads}$ represents the number of threads of the computer. In each group, all threads simultaneously perform exploration at the same time. The number of explorations is a multiple of computer threads; however, the time is the same as that of a single exploration. The reduction in time efficiency is directly proportional to the number of threads in the computer.



\section*{Acknowledgements} This paper is supported by the National Key R\&D Program of China [grant number 2021YFA1000403]; the National Natural Science Foundation of China [grant number 11991022];  the Strategic Priority Research Program of Chinese Academy of Sciences [grant number XDA27000000]; and the Fundamental Research Funds for the Central Universities.


\bibliography{mybibfile}

\end{document}